\documentclass[final,acmtoms]{acmsmallNEW}

\usepackage[ruled]{algorithm2eNew}
\renewcommand{\algorithmcfname}{ALGORITHM}
\SetAlFnt{\small}
\SetAlCapFnt{\small}
\SetAlCapNameFnt{\small}
\SetAlCapHSkip{0pt}
\IncMargin{-\parindent}

\usepackage{algorithmic}

\newcommand{\QIBSH}{\texttt{QIBSH}}
\newcommand{\QIBSHPP}{\texttt{QIBSH++}}
\newcommand{\CPP}{\texttt{C++}}

\newcommand{\MATLAB}{\texttt{Matlab}}
\newcommand{\PACKAGE}[1]{\texttt{#1}}
\newcommand{\LANGUAGE}[1]{\texttt{#1}}
\newcommand{\CLASS}[1]{\textsc{#1}}

\acmVolume{}
\acmNumber{}
\acmArticle{arXiv~}
\acmYear{}
\acmMonth{}

\usepackage{setspace}
\usepackage[noautolanguage]{numprint}
\npstyleenglish

\RequirePackage{amsfonts,amssymb,amsbsy,amsmath}
\usepackage{listings}

\usepackage{multirow}
\usepackage{epsfig,longtable}
\usepackage{soul}
\usepackage{bm}

\usepackage[hyphens]{url}

\usepackage{subfigure}

\usepackage{listings}
\usepackage{color} 
\definecolor{mygreen}{RGB}{14,102,0} 
\definecolor{mylilas}{RGB}{170,55,241}

\newtheorem{rem}{\sc Remark}[section]

\def \c {\bm{c}}

\def \v {\bm{v}}

\def \RR {{\mathbb{R}}}   

\makeatletter
\def\maketag@@@#1{\hbox{\m@th\normalfont\normalsize#1}}
\makeatother

\begin{document}

\lstset{language=Matlab,%
    breaklines=true,%
    morekeywords={matlab2tikz},
    keywordstyle=\color{blue},%
    morekeywords=[2]{1}, keywordstyle=[2]{\color{black}},
    identifierstyle=\color{black},%
    stringstyle=\color{mylilas},
    commentstyle=\color{mygreen},%
    showstringspaces=false,
    numbers=left,%
    numberstyle={\tiny \color{black}},
    numbersep=9pt, 
    emph=[1]{for,end,break},emphstyle=[1]\color{red}, 
}

\markboth
{E. Bertolazzi and A. Falini and F. Mazzia}
{The Object Oriented \CPP{}
 library \QIBSHPP{} for Hermite spline Quasi Interpolation}

\title{The Object Oriented \CPP{}
 library \QIBSHPP{} for Hermite spline Quasi Interpolation}
\author{
  ENRICO BERTOLAZZI
  \affil{%
    Dipartimento di Ingegneria Industriale             
    Universit\`a degli Studi di Trento.
  }   
  ANTONELLA FALINI
  \affil{%
    Dipartimento di Informatica,
    Universit\`a degli Studi di Bari Aldo Moro, Italy.
  }
  FRANCESCA MAZZIA
  \affil{%
    Dipartimento di Informatica,
    Universit\`a degli Studi di Bari Aldo Moro, Italy.
  }
}

\begin{abstract}
  The library \QIBSHPP{} is a \CPP{} object oriented library for the solution of Quasi Interpolation problems.
  The library is based on a Hermite Quasi Interpolating operator, which  
  was derived as continuous  extensions of linear multistep methods 
  applied for the numerical solution of Boundary Value Problems 
  for Ordinary Differential Equations. 
  The library includes the possibility to use Hermite data or to apply 
  a finite difference scheme for derivative approximations,
  when derivative values are not directly available.
  The generalization of the quasi interpolation procedure
  to surfaces and volumes approximation by means of a tensor
  product technique is also implemented. 
  The  method has been also generalized for one dimensional
  vectorial data, periodic data, and for two dimensional data
  in cylindrical coordinates, periodic with respect to the angular argument.
  Numerical tests show that the library could be used efficiently 
  in many  practical problems.
\end{abstract}





\begin{bottomstuff}
Author's addresses:

E. Bertolazzi, 
Dipartimento di Ingegneria Industriale, Universit\`a degli Studi di Trento,
Via Sommarive 9, Trento (ITaly)
email:enrico.bertolazzi@unitn.it

and

A. Falini,
Dipartimento di Informatica,
Universit\`a degli Studi di Bari Aldo Moro, 
Via Orabona  4,  70125 Bari (Italy).
email:antonella.falini@uniba.it

and

F. Mazzia, Dipartimento di Informatica, Universit\`a degli Studi di Bari Aldo Moro, Via Orabona 4, 70125 Bari (Italy).
email:francesca.mazzia@uniba.it
\end{bottomstuff}

\maketitle

\section{Introduction}

There has been a lot of study in constructing good software for
interpolation and  data fitting using spline.
The first package~\PACKAGE{pppack} of de Boor, was available 
in~\PACKAGE{Netlib} from 1992 but the first release was dated 
1971~\cite{DeB2,DeB1}.
The~\MATLAB{} package for spline interpolation and fitting data 
is based on the de Boor subroutines.
Later on in 1973, the algorithm numbered 461, was published on 
ACM Transaction of mathematical software and it was related to
he computation of a cubic spline approximation to the solution
of a linear second order boundary value ordinary differential
equations \cite{burkowski1973algorithm}.
The same journal published in 2016 a B-spline Adaptive Collocation
software for PDEs with Interpolation-Based Spatial 
Error Control \cite{pew2016algorithm}.
Both algorithms use spline functions for the solution of differential
problems with collocation.
In 1993 the \PACKAGE{tspack} package for tension spline curve-fitting 
package~\cite{renka1993algorithm} and  in 2009 its extension for curve
design and data fitting~\cite{renka2009algorithm} have been published.
Nowadays, a lot of wrappers or re-implementations of the \PACKAGE{pppack}
library are available in different languages like~\LANGUAGE{C}, 
\CPP{}, \LANGUAGE{Python}, \MATLAB{}.
Other functions for scattered data are available.
We recall the~\LANGUAGE{Fortran} package~\PACKAGE{fitpack}
of Paul Dierks available in~\PACKAGE{Netlib}~\cite{website:Dierckx}
and the~\LANGUAGE{C} package~\PACKAGE{TSFIT} for two-stage scattered
data fitting \cite{website:Davydov}; the \CPP{} library 
\PACKAGE{G+SMO}~\cite{juttler2014geometry+,mantzaflaris2019overview}
and the \PACKAGE{GeoPDEs} package~\cite{de2011geopdes,vazquez2016new},
both for iso-geometric analysis.
Other libraries and software's available for spline fitting 
and geometric spline constructions
are~\cite{website:Elber,website:Schumaker,SPLINTER,walker2019arbtools} .

The library we present here is  based on the so called 
Hermite BS quasi-interpolant  (BSH QI in short) introduced in~\cite{Mazzia1}, derived from a 
class of linear multistep boundary value methods based on spline
collocation~\cite{Mazzia.Sestini.Trigiante.06.SIAM}.

Univariate spline Quasi Interpolants (QIs) are operators
for function approximations with the following form:
\begin{equation}
  Q_d\left(f\right)=\sum_{j\in J}{\mu_j(f)B_j},
  \label{eq:splineQI}
\end{equation} 
\noindent where $\left\{B_j,j\in J\right\}$ is the B-spline basis
of a given  degree $d$, and $\mu_j(f)$ are local linear functionals.
One of the main properties of QIs is that the coefficients
$\mu_j(f)$ depend locally on the data, making them competitive 
with respect to global approximation methods. 

The  library~\QIBSHPP{} is an object oriented extension of 
the~\LANGUAGE{C} library~\QIBSH{} presented in~\cite{rapportoQIBSH},
\cite{IurinoTesi}
and includes all the procedures for the~\PACKAGE{BSH} Quasi Interpolation
scheme, a generalization of the former to be used when derivative 
values are not available, and an extension of the \PACKAGE{BSH} QI
operator to bivariate  and trivariate functions, which uses a 
suitable tensor product technique. 
Moreover,~\MATLAB{} and~\LANGUAGE{Octave} interfaces have been
implemented for all the objects, in order to make them available
in this well known numerical computing environment.
This is an important feature, since some of the procedures 
in~\QIBSHPP{} will be  also part of the~\MATLAB{} code~\PACKAGE{TOM}
for the numerical solution of Boundary Value Problem for 
Ordinary Differential Equations~\cite{Mazzia4,Mazzia5,Mazzia3}.

The aim of this library is to make available to a wider audience
quasi-interpolation procedures that could be useful when interpolation
is not necessary and the error in the data is negligible.
In many applications, moreover, the first derivative is a known data
and so Hermite quasi-interpolation could give more accurate results
than standard quasi-interpolation.
We experienced a lack of general purpose codes based on high order
quasi-interpolation, especially  for two  and three dimensional data
and in many applications where is required in output continuity
for higher derivatives and for which codes that are based on 
radial basis functions, or bi-variate splines are not suited. 

In Section~\ref{BSH1D} we give a brief description of 
the~\PACKAGE{BSH} QI in one dimension, introducing also 
the approximated~\PACKAGE{BSH}, where derivative values
are not directly used in the operator, but derived using 
suitable finite difference schemes. 
In Section~\ref{sezione_prodotto_tensore} the~\PACKAGE{BSH} 
Quasi Interpolant is extended to the approximation of tensor product
 surfaces and volumes.
In Section~\ref{descrizione_libreria} we describe the implementation
details of the algorithm.
Finally, in Section~\ref{test} we provide some numerical
examples giving an idea of the performance of the~\QIBSHPP{} library. 
The behavior of~\QIBSHPP{} is compared to the QI method \cite{SAB2},
and QI linear, both implemented by the authors and to the spline
interpolation routines from~\MATLAB{}, 
on standard test functions from the literature. 
We also show how to improve the time efficiency of the~\PACKAGE{TOM}
code for BVP problems using the \QIBSHPP{} library.
Moreover, a surface parameterization with high smoothness for complex geometries is presented in subsection \ref{sub:surface}.
We conclude the work showing that the~\QIBSHPP{} library can be applied
for the solution of  two real data problems: a continuous digital elevation model and a biomedical
application.

\section{BSH Quasi Interpolants in One Dimension}\label{BSH1D}

Differential quasi interpolants (DQI)~\cite{DeB4,DeB1} are linear 
approximating operators where the coefficients of the approximating 
splines are computed by linear combinations or averages of derivative
values of $f$, a continuous function defined on an interval $[a,b]$. 
The idea of applying a Hermite Quasi Interpolating technique to 
our problem comes from a different area, since BS methods are a 
class of Boundary Value Methods for ODEs~\cite{brugnano,Mazzia4}. 
Using this class of BS methods it is possible to determine a 
spline $s=\sum_{i\in I}c_iB_i$ on the  mesh 
 defined by the knot vector
$\pi=[x_0,\ldots,x_N]$, where $a=x_0<x_1<\ldots<x_N=b$, satisfying the Hermite interpolation 
conditions $s(x_i)=f_i$, $s'(x_i)=f'_i$ for all $i=0,\ldots,N$, 
where  $f_0,\ldots,f_N$ and $f'_0,\ldots,f'_N$ are respectively
the values of the function $f$ and of the first derivative 
$f^\prime$ and both, $f_i,f'_i\in\RR^{p}\; \forall i$, with $p>1$ for the multidimensional case.

Here, the set $\left\{B_i:\;i\in I\right\}$ is the B-spline basis for 
the space  $S_{d,\pi}$  of  $d$-degree splines on the knots $\pi$. 
The BS Hermite Quasi Interpolation scheme approximates a function 
$f$ on an interval $\left[a,b\right]$ starting from its values,
and from those of the first derivative on $N+1$ mesh points
$\pi$.
We want an approximating function in the space 
$S_{d,\pi}$ of the splines of degree $d$ with knots $\pi$. 
Usually, we work with an extended knot set considering a total of $2d$
additional boundary knots.
The new knot set is then defined as   
$\boldsymbol{\tau}=\{\tau_i\}_{i=1}^{N_\tau}=\left\{x_{-d},\ldots,x_{-1},x_0,\ldots,x_N,x_{N+1},\ldots,x_{N+d}\right\}$,
where $N_\tau=N+2d+1$.
The auxiliary boundary knots are commonly chosen equal to the ending
points of the interval. The Quasi Interpolating spline is then:
\begin{equation}
  Q^{\left(BS\right)}_d\left(f\right)
  =\sum_{j=-d}^{N-1}{\mu^{\left(BS\right)}_j(f)B_j},
  \label{eq:spline1D}
\end{equation} 
\noindent with the coefficients $\mu_j^{(BS)}(f)$ expressed by 
\begin{equation}
  \mu^{(BS)}_j(f)=
  \displaystyle\sum_{i=1}^{d}
  \begin{cases}
  \hat{\alpha_i}^{(-1,j+d+1)}f_{i-1}-h_{k_1}
  \hat{\beta_i}^{(-1,j+d+1)}f'_{i-1} 
  & j=-d,\ldots,-2, 
  \\[1em]
  \hat{\alpha_i}^{(j,d)}f_{i+j}-h_{j+k_1+1}
  \hat{\beta_i}^{(j,d)}f'_{i+j}
  & j=-1,\ldots,\tilde{N},
  \\[1em]
  \hat{\alpha_i}^{(\tilde{N},j+d-\tilde{N})}f_{\tilde{N}+i}
  -h_{N-k_2}\hat{\beta_i}^{(\tilde{N},j+d-\tilde{N})}f'_{\tilde{N}+i}
  & j=\tilde{N}+1,\ldots,N-1,
\end{cases}
\label{eq:mujlocal}
\end{equation}

\noindent  where $\tilde{N}=N-d$, $\bm{\hat{\alpha}}^{(j,r)}=\left(\hat{\alpha}^{(j,r)}_{1},\ldots,\hat{\alpha}^{(j,r)}_{d}\right)^{T}$ and $\bm{\hat{\beta}}^{(j,r)}=\left(\hat{\beta}^{(j,r)}_{1},\ldots,\hat{\beta}^{(j,r)}_{d}\right)^{T}$ are solutions of local linear systems  of size $2d\times2d$, whose coefficient matrix depends on the values of the B-splines $B_i$ (see \cite{Mazzia4}). 
The functionals $\mu_{j}^{(BS)}(f)$ in \eqref{eq:spline1D} depend 
locally on the values $(f_i,f'_i),\,i=0,\ldots,N$,
where $f_i$ and $f'_i$ denote the exact function and 
first derivative values, respectively.
Following the notation used in~\cite{Mazzia2},
we write the coefficients in the form:
\begin{equation}\label{eq:mujlocalmatrixform}
  \bm{\mu}^{(BS)}=(\hat{A}\otimes I_{p})\mathbf{f}-(\hat{H}\hat{B}\otimes  I_{p})\mathbf{f}',
\end{equation}
where $\hat{A}$ and $\hat{B}$ are banded matrices in  $\RR^{N+d\times N+1}$
containing the local coefficients $\bm{\hat{\alpha}}^{(j,r)}$ and $\bm{\hat{\beta}}^{(j,r)}$, while $\hat{H}=\mathrm{diag}(\hat{h}_1,\ldots,\hat{h}_{N+d})$, with

\begin{equation}\label{eq:hi}
  \hat{h}_i=
  \begin{cases}
    {h}_{k_1},\quad&\text{if}~~i\leq d\\
    {h}_{k_1+i-d},\quad&\text{if}~~d+1\leq i \leq N\\
    {h}_{N-k_2},\quad&\text{if}~~i\geq N+1,
  \end{cases}
\end{equation}
and  $\mathbf{f}=\left(f_{0}^{\top},\ldots,f_{N}^{\top}\right)^{\top}$,
$\mathbf{f}'=({f'_{0}}^{\top},\ldots,{f'_{N}}^{\top})^{\top}$, 
and $I_{p}$ is the identity operator of dimension $p$.

The coefficient determination of the Quasi Interpolant
$Q^{(BS)}_d(f)$ in the B-spline representation, requires 
the solution of local linear systems for $j=-1,\:r=1,\ldots,d$,
for $0\leq j\leq N-d-1,\:r=d$, and for $j=N-d,\:r=d,\ldots,2d-1$.
The computational cost is the one for solving $N+d$ linear 
systems of dimension $2d\times2d$. 
They are solved by the efficient and stable algorithm 
presented in~\cite{Mazzia5}.
When low degree polynomials are used, we can give explicit expressions 
of the coefficient vectors.
Note that for a uniform knot mesh $\pi$, the inner coefficient
vectors $\bm{\hat{\alpha}}^{(j,d)}$ and $\bm{\hat{\beta}}^{(j,d)}$,
for $j=-1,\ldots,N-d$ do not depend on $j$,
so the expressions of coefficients $\mu^{(BS)}_{j}(f)$ 
may be derived  beforehand.

The Quasi Interpolation procedure described so far is of Hermite type,
since it depends on the function and its first derivative values.
Often in applications, we do not have such information, and only
approximate values of $f$ and $f'$ are available.
In some cases, we may be given only the values of the function
at mesh knots, and in order to construct the QI we must use
approximate values for the first derivatives.
For this reason, we combine the Quasi Interpolation scheme with a 
symmetric finite difference scheme approximating the derivatives
of $f$ at the mesh points.
In order to distinguish it from the original one, 
we refer to the original BSH as $Q_d^{(BS)}$,
and to the one using approximate values of $f'$, as $Q^{(BSa)}_d$. 

Indeed, we can approximate the first derivative values of a 
sufficiently smooth function $f$ on a grid $x_0<\ldots<x_N$,
using the $l$-step finite difference scheme used in~\cite{Mazzia2}.
So we have  the following scheme:
\begin{equation}\label{eq:finitedifferencescheme}
  \begin{cases}
    \displaystyle\sum_{i=0}^{l}\gamma_{i}^{(n)}f_i=f'(x_n)+\mathcal{O}(h^l),
    \quad & n=0,\ldots,l_1-1,\\
    \displaystyle\sum_{i=0}^{l}\gamma_{i}^{(n)}f_{n-l_1+i}=f'(x_n)+\mathcal{O}(h^l),
    \quad & n=l_1,\ldots,N-l_2,\\
    \displaystyle\sum_{i=0}^{l}\gamma_{i}^{(n)}f_{N-l+i}=f'(x_n)+\mathcal{O}(h^l),
    \quad & n=N-l_2+1,\ldots,N,
  \end{cases}
\end{equation}
where $l_1=\lfloor{l/2}\rfloor$, $l_2=l-l_1$, $h_i=x_i-x_{i-1}$, 
$i=1,\ldots,N$, and $h=\max_{1\leq i\leq N}h_i$.
The coefficients $\gamma_i^{(n)}$ are computed imposing that the 
local truncation error of the resulting method is $\mathcal{O}(h^l)$.
The derivative approximation scheme is 
modified in order to have a symmetric global approximation 
when the mesh is symmetric. It can be  written as 
\begin{equation}\label{eq:ordineGBDF}
  (\Gamma^{(l)}\otimes I_{p})
  \begin{pmatrix}
    f(x_0)\\
    \vdots\\
    f(x_N)
  \end{pmatrix}=
  \begin{pmatrix}
    f'(x_0)\\
    \vdots\\
    f'(x_N)
  \end{pmatrix}+
  \mathcal{O}(h^l),
\end{equation}
where the banded matrix $\Gamma^{(l)}\in\RR^{(N+1)\times (N+1)}$ 
has the following structure: 
\begin{equation}\label{eq:Gamma}
  \Gamma^{(l)}=\begin{pmatrix}
  \gamma_0^{(0)} & \cdots & \gamma_l^{(0)}&0&\cdots&\cdots&0      \\
  \vdots         &        & \vdots&\vdots&  &      &      \vdots \\
  \gamma_0^{(l_1-1)}& \cdots &\gamma_l^{(l_1-1)}&0&&&\vdots\\
  0&\ddots & &\ddots&\ddots&                        &\vdots     \\
  \vdots&\ddots&\gamma_0^{(j)}& \cdots &\gamma_l^{(j)}&\ddots  &\vdots    \\
  \vdots&&\ddots&\ddots&&\ddots&0                      \\
  \vdots&&&0&\gamma_0^{(N-l_2+1)} & \cdots & \gamma_l^{(N-l_2+1)}      \\
  \vdots&&&\vdots&\vdots         &        & \vdots              \\
  0&\cdots&\cdots&0&\gamma_0^{(N)}& \cdots &\gamma_l^{(N)}
  \end{pmatrix}.
\end{equation}

The first derivative approximation on the mesh points is 
$\mathbf{\tilde{f}'}=(\Gamma^{(l)}\otimes I_{p})\mathbf{f}$.
Combining this scheme with the BSH formula for the
QI coefficients we get the Quasi Interpolant:
\begin{equation}\label{eq:QIBSHa1D}
  Q^{\left(BSa\right)}_d\left(f\right)
  =\sum_{j=-d}^{N-1}{\tilde{\mu}^{\left(BSa\right)}_j(f)B_j},
\end{equation}
It can be proven that the following error bound holds:
\begin{equation}
  \|f-Q^{(BSa)}_d(f)\|_{\infty}\leq Lh^{d+1}\|D^{d+1}f\|_{\infty}
  +\tilde{L}\,h^{l+1},
\end{equation}
where $L$ and $\tilde{L}$ depends on the spline degree $d$,
on the smoothness of $f$ and on the geometric properties,
like quasi-uniformity, of the underlying knot mesh.
For more details we refer to~\cite{Falini.Mazzia.Sestini.2022}.
Note that a similar error bound holds also for the operator $Q_d^{(BS)}$,
see~\cite{Mazzia.Sestini.09Bit} for the details.

\section{BSH Quasi Interpolants and Tensor Product}
\label{sezione_prodotto_tensore}

The general tensor product framework can be used with the BSH QI
from the previous section.
Applying the one dimensional operators  either first along 
the x-direction and then the y-direction, or vice-versa,
we compute the approximating surface.
Suppose we are given a function $f:\RR^2\rightarrow\RR$ 
$f=f(x,y)$ defined  on a rectangular planar domain 
$R=\left[a,b\right]\times\left[c,d\right]$.
We choose the spline degrees to be $d_x$ and $d_y$ and select 
$m_1=N+1$ and $m_2=M+1$ knots respectively on $x$ and $y$ directions.
Consider the partitions with the additional boundary knots defined as:
\begin{equation}
  \begin{array}{l}
	 \tau_x=x_{-d_x}\leq\cdots\leq x_{0}=a<\cdots x_i\cdots  <b=x_{N}\leq x_{N+1}\leq\cdots\leq x_{N+d_x},\\[1em]
	 \tau_y=y_{-d_y}\leq\cdots\leq y_{0}=c<\cdots y_j\cdots <    d=y_{M}\leq y_{M+1}\leq\cdots\leq y_{M+d_y}. 
  \end{array}
\end{equation}
We denote the extended knot vectors as $\bm{\tau}_x$ and $\bm{\tau}_y$, 
with sizes  respectively $N_{\tau_x}$ and $N_{\tau_y}$.
The dimension of the spline spaces are $n_1=N+d_x$ and $n_2=M+d_y$
respectively for each knot partition.
Any spline $s$ in the tensor product space is written in the form:
\begin{equation}\label{eq:spazioS1S2}
  s(x,y)=\sum^{n_1}_{p=1} \sum^{n_2}_{q=1} c_{pq}\phi_{p}(x)\psi_{q}(y)\quad\forall (x,y) \in \RR^2,
\end{equation}
where $\phi_{p}$ and $\psi_{q}$ are the B-splines elements.
Also, the exact values of the function $f$ and of its partial derivatives
$f_x$,$f_y$ and $f_{xy}$ are assigned on the grid points given in matrix
form respectively as $\mathbf{F}=(f)_{i,j}$, $\mathbf{F}_x=(f_x)_{i,j}$,
$\mathbf{F}_y=(f_y)_{i,j}$, and 
$\mathbf{F}_{xy}=(f_{xy})_{i,j}\in\RR^{m_1 \times m_2}$, where
$f_{i,j}=f(x_i,y_j)$ for all indices.
The aim is to find a function $g=g(x,y)$ in the tensor product space~\eqref{eq:spazioS1S2} aproximating $f$,
finding the coefficient matrix $\mathbf{C}=(c_{pq})\in \RR^{n_1\times n_2}$.
The tensor product technique applied to the BSH follows a very simple idea: 
the grid is divided into knot lines according to its knot partitions,
and a specific number of one dimensional QI problems is solved 
along the two directions. 

\noindent The coefficients of the spline function $s(x,y)$ are computed 
using the QI BSH tensor product algorithm, which is summarized 
by the following steps.

\begin{longenum}

\item For each $j=1,...,m_2$, pick up the function  and the $x$ partial
      derivatives of the $j$-th line in the grid, and use
      them applying the one dimensional BSH $Q_{d_x}^{(BS)}$.
      Solve $m_2$ one dimensional QI problems, storing the 
      resulting coefficients into the matrix 
      $\mathbf{D}=(d_{pj})\in\RR^{n_1\times m_2}$.

      In matrix formulation we have
      \begin{equation}\label{eq:D}
         \mathbf{D}=\hat{A}_x\mathbf{F}-\hat{H}_x\hat{B}_x\mathbf{F}_x,
      \end{equation}
      where the subscript in the matrices $\hat{A}_x$, $\hat{B}_x$ 
      and $\hat{H}_x$,  specifies the equation~\eqref{eq:mujlocalmatrixform}
      for the local coefficients on the $x$-direction with $p=1$. 

\item For each $j=1,...,m_2$, pick up the $y$ and the mixed $xy$ 
      partial derivatives of the function $f$ on the $j$-th 
      line in the grid, and apply the one dimensional operator
      $Q_{d_x}^{(BS)}$.
      As in the previous step, $m_2$ one dimensional QI problems
      are solved, saving the coefficients into
      $\mathbf{D'}=(d'_{pj})\in\RR^{n_1\times m_2}$.
      In matrix formulation:
      \begin{equation}\label{eq:Duno}
         \mathbf{D'}=\hat{A}_x\mathbf{F}_y-\hat{H}_x\hat{B}_x\mathbf{F}_{xy}.
      \end{equation}

\item Now, simply switch the direction, and consider the quasi
      interpolant operator $Q_{d_y}^{(BS)}$ along the $y$ direction.
      For each $p=1,\ldots,n_1$ pick up the values in the $p$-th rows
      of the matrices $\mathbf{D}=(d_{p,j})$ and $\mathbf{D'}=(d'_{p,j})$,
      and apply to these values the one dimensional BSH $Q_{d_y}^{(BS)}$.
      So, after solving $n_1$ one dimensional QI problems in the 
      $y$ direction, the matrix $\mathbf{C}=(c_{p,q})$ containing
      the coefficients of the tensor product form~\eqref{eq:spazioS1S2} is built. 
      Note that in this last step we are using the BSH quasi interpolating
      scheme $Q_{d_y}^{(BS)}$ on the $y$-direction, replacing the function
      and derivative values by the elements from the matrices
      $\mathbf{D}$ and $\mathbf{D}'$, namely in matrix form: 
      \begin{equation}
         \mathbf{C}=
         \mathbf{D}\hat{A}^{T}_y-\mathbf{D'}(\hat{H}_y\hat{B}_y)^{T}
         =\left(\hat{A}_x\mathbf{F}-\hat{H}_x\hat{B}_x\mathbf{F}_x\right)
         \hat{A}^{T}_y-\left(\hat{A}_x\mathbf{F}_y
         -\hat{H}_x\hat{B}_x\mathbf{F}_{xy}\right)
         \left(\hat{H}_y\hat{B}_y\right)^{T}.
      \end{equation}
\end{longenum}
When only function values are available, approximations
for partial derivatives in both directions are computed
using the scheme~\eqref{eq:ordineGBDF}, that is 
\begin{equation}\label{eq:FxFyFxy}
  \mathbf{F}_x\approx \Gamma^{(l_x)}_{x} \mathbf{F}, \quad \mathbf{F}_y\approx\mathbf{F}(\Gamma^{(l_y)}_{y})^{T}, \quad \mathbf{F}_{xy}\approx\mathbf{F}_x(\Gamma^{(l_y)}_y)^{T}\approx\Gamma^{(l_x)}_x
  \mathbf{F}_y  \approx\Gamma^{(l)}_x\mathbf{F}(\Gamma^{(l_y)}_y)^{T}.
\end{equation}
Matrix expressions~\eqref{eq:D} and~\eqref{eq:Duno} from the 
tensor product scheme are now
\begin{equation}\label{eq:Da}
  \mathbf{D}^{(a)}
  =\hat{A}_x\mathbf{F}-\hat{H}_x\hat{B}_x\Gamma^{(l_x)}_{x} \mathbf{F},
  \quad\text{ and }\quad\mathbf{D'}^{(a)}
  =\hat{A}_x \mathbf{F}(\Gamma^{(l_y)}_{y})^{T}-\hat{H}_x\hat{B}_x
  \Gamma^{(l_y)}_x\mathbf{F}(\Gamma^{(l_y)}_y)^{T},
\end{equation}
and they satisfy
\begin{equation}\label{eq:Dunoa}
  \mathbf{D'}^{(a)}=\mathbf{D}^{(a)}(\Gamma^{(l_y)}_y)^{T}.
\end{equation}
This proves it is theoretically equivalent, to compute the
approximation of partial derivatives starting from $\mathbf{F}$,
and then apply the tensor product of  BSH QI, or apply the
finite difference scheme to the elements in  $\mathbf{D}^{(a)}$.
The latter way has been used in our library implementation, 
since it reduces the computational cost.
The coefficients matrix of the representation~\eqref{eq:spazioS1S2}
is obtained in this case as
\begin{equation}
  \mathbf{C}^{(a)}=(\hat{A}_x-\hat{H}_x\hat{B}_x\Gamma^{(l_x)}_{x})
  \mathbf{F}\hat{A}_y^{T}-(\hat{A}_x-\hat{H}_x\hat{B}_x
  \Gamma^{(l)}_{x})\mathbf{F}(\Gamma^{(l_y)}_y)^{T}(\hat{H}_y\hat{B}_y)^{T},
\end{equation}
or, in a more compact form:
\begin{equation}\label{eq:coeff2D}
   \mathbf{C}^{(a)}=(\hat{A}_x-\hat{H}_x\hat{B}_x
   \Gamma^{(l_x)}_{x})\mathbf{F}(\hat{A}_y-\hat{H}_y\hat{B}_y
   \Gamma^{(l_y)}_{y})^{T}.
\end{equation}
In the following, we address the extension of the operator
$Q_d^{(BSa)}$ to the approximation of a three valued 
function $f=f(x,y,z)$, with  $f:\RR^3\rightarrow\RR$  
defined on a domain
$R=\left[a_x,b_x\right]\times\left[a_y,b_y\right]\times\left[a_z,b_z\right]$.

Consider the partitions with the additional boundary knots defined as:
\begin{equation}
  \begin{array}{l}
	 \tau_x=x_{-d_x}\leq\ldots\leq x_{0}=a_x<\ldots x_i\ldots  <b_x=x_{N_x}\leq x_{N_x+1}\leq\ldots\leq x_{N_x+d_x},\\[1em]
	 \tau_y=y_{-d_y}\leq\ldots\leq y_{0}=a_y<\ldots y_j\ldots <    b_y=y_{N_y}\leq y_{N_y+1}\leq\ldots\leq y_{N_y+d_y}. \\[1em]
	 	 \tau_z=z_{-d_z}\leq\ldots\leq z_{0}=a_z<\ldots z_j\ldots <    b_z=z_{N_z}\leq z_{N_z+1}\leq\ldots\leq z_{N_z+d_z}. 
  \end{array}
\end{equation}
To explain better the application of the quasi-interpolant on
a multi dimensional space it is convenient to express the
data using a three dimensional tensor 
\[
  \mathcal{F} \in \RR^{Nx \times Ny \times Nz}, \qquad \mathit{f}_{ijk}=f(x_i,y_j,z_k) 
\]
and the coefficients of the  spline using  $n$-mode product
between a tensor $\mathcal{X} \in \RR^{I_1 \times I_2 \times I_3}$
and a matrix $A \in \RR^{n \times I_j}$ as defined in~\cite{kolda}:
\[
  \mathcal{Y} = \mathcal{X} \times_j A \Leftrightarrow  Y_j = A X_j 
\]
where $Y_j, X_j$ are the mode-$j$ unfolding of ${\cal Y}$  and ${\cal X}$.

Here we report only the version of the quasi-interpolant using
approximate derivative, because it is less expensive from a 
computational point of view.
The coefficient of the approximated  BSH QI version can been 
in fact calculated, like for the bi-dimensional case,
by applying the one dimensional quasi-interpolant 
along the three directions. We obtain that:
\begin{equation}\label{eq:C3}
  \mathcal{C}^{(a)}=\mathcal{F} 
  \times_1 (\hat{A}_x-\hat{H}_x\hat{B}_x\Gamma^{(l_x)}_{x})
  \times_2 (\hat{A}_y-\hat{H}_y\hat{B}_y\Gamma^{(l_y)}_{y})^{T}
  \times_3 (\hat{A}_z-\hat{H}_z\hat{B}_z\Gamma^{(l_z)}_{z})^{T}.
\end{equation}
Additional details can be found in \cite{Falini.Mazzia.Sestini.2022}.

\begin{rem}
  Note that by using the $n$-mode product we can express
  the coefficient matrix of equation \eqref{eq:coeff2D} as,
  \begin{align*}
    C^{(a)} = {\cal F} \times_1 (\hat{A}_x-\hat{H}_x\hat{B}_x\Gamma^{(l_x)}_{x})
    \times_2 (\hat{A}_y-\hat{H}_y\hat{B}_y\Gamma^{(l_y)}_{y})^{T}.
  \end{align*}
  With this operation also the generalization to any
  dimension can be easily derived.
\end{rem}

\section{Implementation Details}\label{descrizione_libreria}

The library \QIBSHPP{} is a collection of \CPP{} procedures
managing the Hermite Quasi Interpolation 
for functions of one (scalar and vectorial functions), two  and three variables. 
A \MATLAB{} toolbox is also available, that allow to handle all the procedure in \MATLAB{}.

This library is an improved version of the C
version of \QIBSH{} library described in \cite{IurinoTesi,rapportoQIBSH}.
The new implementation use more efficiently memory
and dynamic memory usage. The C++ classes are mapped 
in \MATLAB{} classes using MEX interfacing mechanism.
In practice each \MATLAB{} class instance store
a pointer to the corresponding \CPP{} class instance
and each method for the \MATLAB{} class call 
a method of the corresponding \CPP{} class.

This approach permits to develop and test algorithm
using \QIBSH{} in \MATLAB{}.
This remapping introduce an overhead that is small 
and acceptable for the proposed applications.
In any case, for best performance it is easy
to translate to \MATLAB{} to \CPP{}.

\subsection{\CPP{} classes}
The library is organized as a set of classes that interact
together for the Quasi interpolant build and evaluation. 
\begin{itemize}
  \item B-spline basis computation 
  \begin{itemize}
    \item \CLASS{Bspline}\\
    This class implement the classical B-spline
    as described in~\cite{DeB1,SC1}.
    Recurrence formula for derivative and integral, standard knot placements, tensor product B-spline evaluation given the support polygon.
  \end{itemize}
  \item Finite Difference
  \begin{itemize}
    \item \CLASS{FiniteDifferenceUniformD1}
    \item \CLASS{FiniteDifferenceD1}\\
    This two classes implement the finite difference 
    approximation of derivative given a list $(x_i,y_i)$ of interpolation points. The class  \CLASS{FiniteDifferenceUniformD1} do the same
    computation more efficiently when coordinated $x_i$ are uniformly distributed. 
  \end{itemize}
  \item Derivative approximation 
  \begin{itemize}
    \item \CLASS{ApproximateDerivative}
    \item \CLASS{ApproximateDerivative2D}\\
    This classes uses \CLASS{FiniteDifferenceUniformD1} and
    \CLASS{FiniteDifferenceD1} to build the finite difference
    approximation of a set of one dimensional $(x_i,y_i)$
    or two dimensional $(x_{ij},y_{ij},z_{ij})$ 
    function sampling. In 2D cases mixed derivatives
    are obtained by applying finite difference in the $y$
    direction to the approximate $x$ derivative.
    The classes manages cyclic data approximation 
    if required.
  \end{itemize}
  \item Quasi Hermite 1D
  \begin{itemize}
    \item \CLASS{QuasiHermite}
    \item \CLASS{QuasiHermiteApprox}\\
    The classes compute the B-spline polygon corresponding 
    to the \QIBSH{} approximation.
    The first class uses points and analytical derivative 
    at the corresponding points.
    The second class approximate the derivative
    using finite difference from class  \CLASS{ApproximateDerivative}.+
  \end{itemize}
  \item Quasi Hermite 2D
  \begin{itemize}
    \item \CLASS{QuasiHermite2D}\\
    This is the base class that compute the polygon for the B-spline that correspond to the \QIBSH{} approximation
    of 2D surface data.
    The 2D points are passed to the class with $x$, $y$ and 
    mixed $xy$ derivatives at the points.
    \item \CLASS{QuasiHermite2Dapprox}\\
    This class is derived from \CLASS{QuasiHermite2D}
    and compute the polygon for the B-spline that correspond to the \QIBSH{} approximation
    of 2D surface data.
    The derivative respect to $x$, $y$ and 
    mixed $xy$ needed for base class are approximated
    using finite difference with class \CLASS{ApproximateDerivative2D}.
    \item \CLASS{QuasiHermite2Dsurface}\\
    This class is derived from \CLASS{QuasiHermite2D}.
    In addition the class store the computed B-spline
    polygon so that can be evaluated at any points
    without the requirement to pass the B-spline
    polygon.
    \item \CLASS{QuasiHermite2DapproxSurface}\\
    This class is derived from \CLASS{QuasiHermite2Dapprox}.
    In addition the class store the computed B-spline
    polygon so that can be evaluated at any points
    without the requirement to pass the B-spline
    polygon.
  \end{itemize}
\end{itemize}

\subsection{\MATLAB{} classes}
\QIBSHPP{} library is connected with \MATLAB{} 
using MEX interface.
The mapping is not one to one but is a little bit 
of higher level.

\begin{itemize}
  \item B-spline basis computation
  \begin{itemize}
    \item \CLASS{Bspline}\\
    This is a one-to-one remap of the corresponding C++ class.
    \item \CLASS{Bspline1D}\\
    This \MATLAB{} class remap the \CPP{} class \CLASS{QuasiHermite} and \CLASS{QuasiHermiteApprox}.
    It stores in the \MATLAB{} class data the B-spline
    polygon and can use analytical derivative or use
    finite difference approximation.

  \end{itemize}
  \item Approximation of derivatives using finite difference
  \begin{itemize}
    \item \CLASS{ApproximateDerivative1D}
    \item \CLASS{ApproximateDerivative2D} \\
    This is a one-to-one remap of the corresponding C++ class.
  \end{itemize}
  \item Spline build using quasi interpolation 
  \begin{itemize}
    \item \CLASS{QIBSH1D}
    \item \CLASS{QIBSH2D} \\
    This are high level remaps of the \CPP{} classes 
    \CLASS{QuasiHermite}, \CLASS{QuasiHermiteApprox},
    \CLASS{QuasiHermite2D}, \CLASS{QuasiHermiteApprox2D}
    with the storage in the \MATLAB{} class of the 
    resulting B-spline polygon.

  \end{itemize}
\end{itemize}

\begin{figure}[!htb]
    \centering
    \includegraphics[width=10cm]{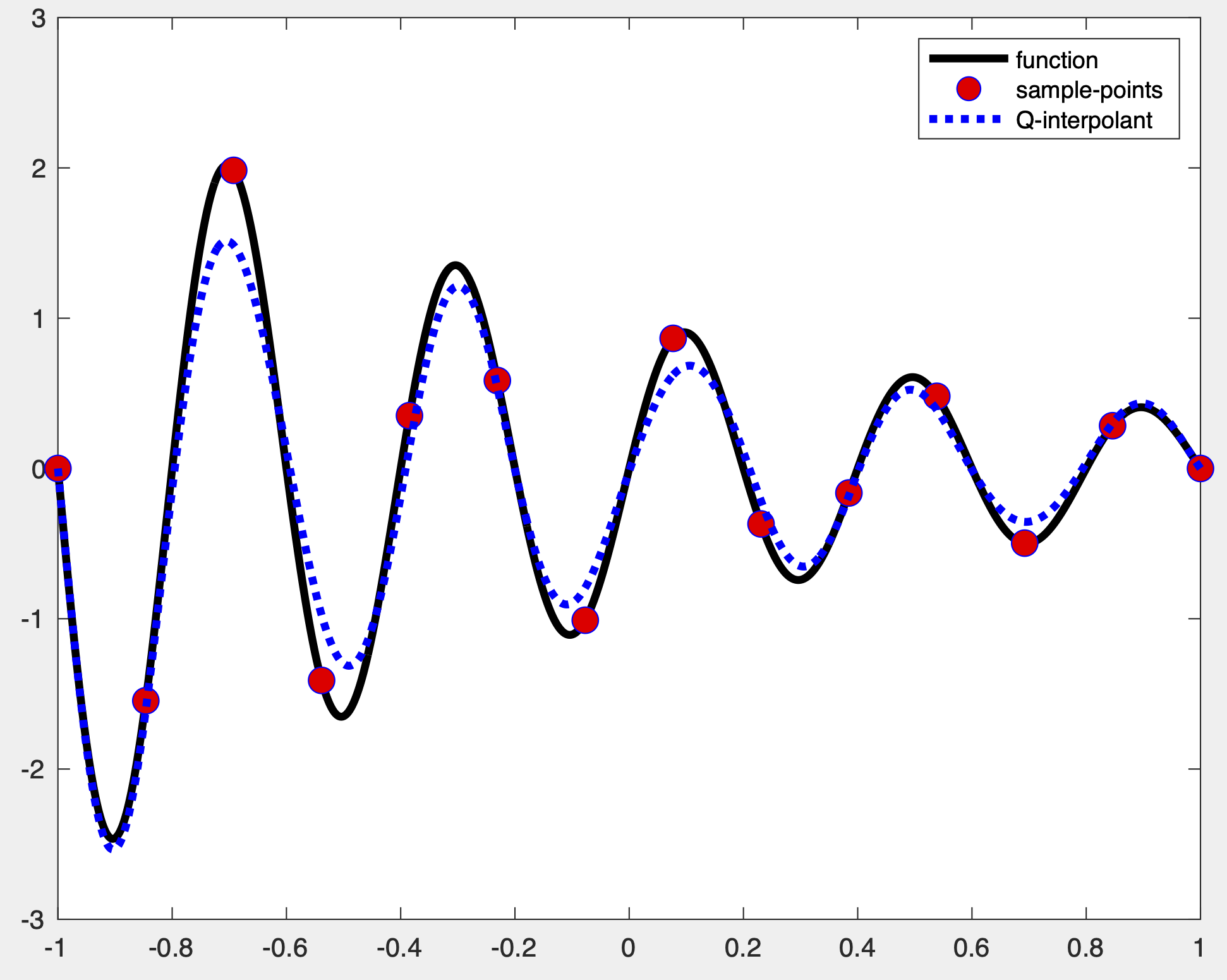}
    \caption{Output of quasi-interpolant for function 
    $\mathrm{e}^{-x}\sin(5\pi x)$ with $14$ equally spaced 
    sample points in the interval $[-1,1]$.}
    \label{fig:qi-example}
\end{figure}

The \MATLAB{} usage of the \QIBSH{} library is
particularly simple:

\begin{lstlisting}[frame=single,language=MATLAB]
% Instantiate in Q a MATLAB class for 1D Quasi Interpolant
Q = QIBSH1D('Qhermite');
% Build a quasi interpolant
% d           = degree of the interpolant
% (x(i),y(i)) = quasi interpolation points
% yprime(i)   = derivative at the quasi interpolation point
% fifth argument if true build a periodic quasi interpolant
Q.build( d, x, y, yprime, false );
\end{lstlisting}

after build is easy to compute points and derivative
on the B-spline:

\begin{lstlisting}[frame=single,language=MATLAB]
% Evaluate y = Q(x)
y = Q.eval(x);
% Evaluate y = Q'(x)
Dy = Q.eval(x,1);
% Evaluate y = Q''(x)
DDy = Q.eval(x,2);
\end{lstlisting}

The interface is very intuitive with few example
it is easy to practice with the library.
Here is a example of quasi-interpolation of
a set of points taken from a sampling of a function:

\begin{lstlisting}[frame=single,language=MATLAB]
% set the function to be approximated
effe  = @(x) exp(-x).*sin(5*pi*x);
effe1 = @(x) ((5*pi*cos(5*pi*x)) - sin(5*pi*x)).*exp(-x);
[a,b] = deal(-1,1); % range of the function
d     = 4;          % degree of quasi interpolant
N     = 14;         % number of sampling points
m     = 1000;       % number of evaluation points for plotting
% plot the function
xx = linspace(a,b,m); yy = effe(xx);
plot( xx, yy, 'LineWidth', 3, 'Color', 'black' );

% Sample function and evaluate f and f' at sample points
x = linspace(a,b,N); y = effe(x); yprime = effe1(x);

% Create a QIBSH object for 1D quasi-interpolation
Q = QIBSH1D('Qhermite');
% build quasi interpolant for (x,y) data, no cyclic data.
Q.build( d, x, y, yprime, false );

% plot the sampling points
hold on;
plot( x, y, 'ob', 'MarkerSize', 10, 'MarkerFaceColor', 'red' );

% plot the approximated interpolant
plot( xx, Q.eval( xx ), ':', 'LineWidth', 3, 'Color', 'blue' );
legend('function','sample-points','Q-interpolant');
\end{lstlisting}
The output of the script is shown in Fig.~\ref{fig:qi-example}.

\subsection{Procedures for Quasi Interpolation of Function of One Variable}\label{procedure_1D}

Using the library we can solve the one dimensional problem of 
Quasi Interpolation. Given the values of a function and of its 
first derivatives at some points, provided a knot set and the
desired spline degree, the user can find the Hermite spline 
quasi interpolating the data. 
The final approximation of the $QI^{(BS)}$ is given in terms of its B-spline coefficients. 
Optionally, function and derivative evaluations can also be returned. When first derivative values are not available,  their approximations are first computed and the $QI^{ (BSa)}$ quasi-interpolant is used.

In addition, specific procedures for the treatment of periodic functions are available. 
These functions specializes the QI procedures for the coefficients and for the finite 
difference scheme, when periodic knots are used.

%
%
\subsection{Procedures for Tensor Product Quasi Interpolation}\label{procedure_2D}
The computation of two dimensional  BSH QI in the tensor product form follows the scheme from Section \ref{sezione_prodotto_tensore}. For data organized on regular grids, the idea is to split the process into one dimensional problems along the axes directions, quasi interpolating the values of the function and of partial derivatives given in matrix form. 
When partial derivatives are not available, 
we combines the tensor product of  BSH  QI with the finite difference scheme \eqref{eq:finitedifferencescheme}. 
The implementation  to reduce the computational costs is based on the computation of  $x$-partial derivatives values 
of the matrix $\mathbf{F}$ and generating the approximations for $\mathbf{F_x}$. Then, one dimensional QI problems along the $x$ axis are solved 
storing the coefficients into an auxiliary matrix $\mathbf{D}^{(a)}$. 
Now, the approximate derivative of  
 $\mathbf{D}^{(a)}$ are computed and stored in   the matrix $\mathbf{D'}^{(a)}$. These two matrices are used for the last step: the solutions of  QI problems along the $y$ axis, having $\mathbf{D}^{(a)}$ and $\mathbf{D'}^{(a)}$ as inputs, returns the spline coefficients. Computational cost is reduced since the aprroximation of the derivatives 
 are related to two function calls rather than  three.
 Moreover, 
for the matrix case, the computation of the coefficients $\gamma_i^{(n)}$ in \eqref{eq:finitedifferencescheme}, and of the matrix  $\Gamma^{(l)}$ in \eqref{eq:ordineGBDF} is done only once.
Finally we remark that also the tensor product case has been generalized
and adapted  for function domains which can easily be described in polar coordinates in the plane.
Hence, 
we can  use the library  to compute the QI or derivative approximations for a periodic 
function $f=f(\rho, \theta)$ of period $T>0$ with respect to its second argument $\theta$.

In the $m$D case the quasi interpolant is implemented only when the derivatives are approximated, using a generalization of the 2D procedure. In this case the computational costs remain related to the dimension $m$. The use of the exact derivatives has a higher computational cost, so it has not be considered.

\begin{figure}[!htb]
    \centering
    \includegraphics[width=6.5cm]{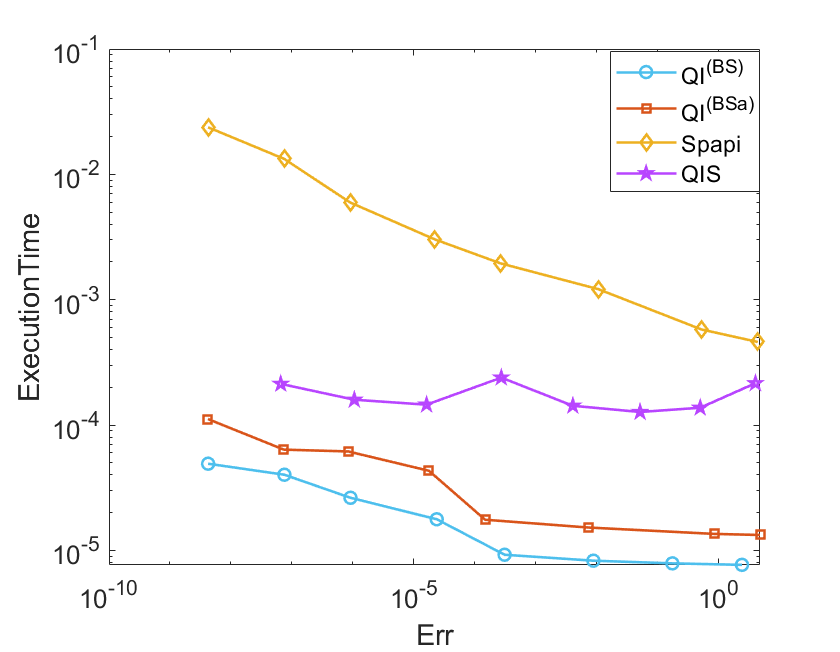}
    \includegraphics[width=6.5cm]{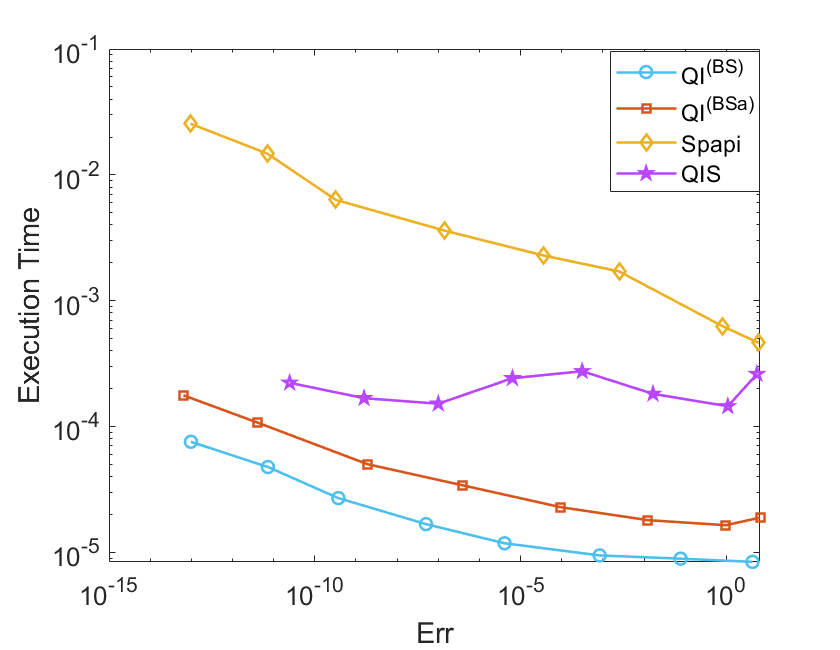}
    \caption{Work precision diagram for the function $f_1$. On the left $d=3$ is used, while on the right $d=5$ is adopted.}
    \label{fig:tempi1D}
\end{figure}

\begin{figure}
    \centering
 
    \includegraphics[height=5.15cm]{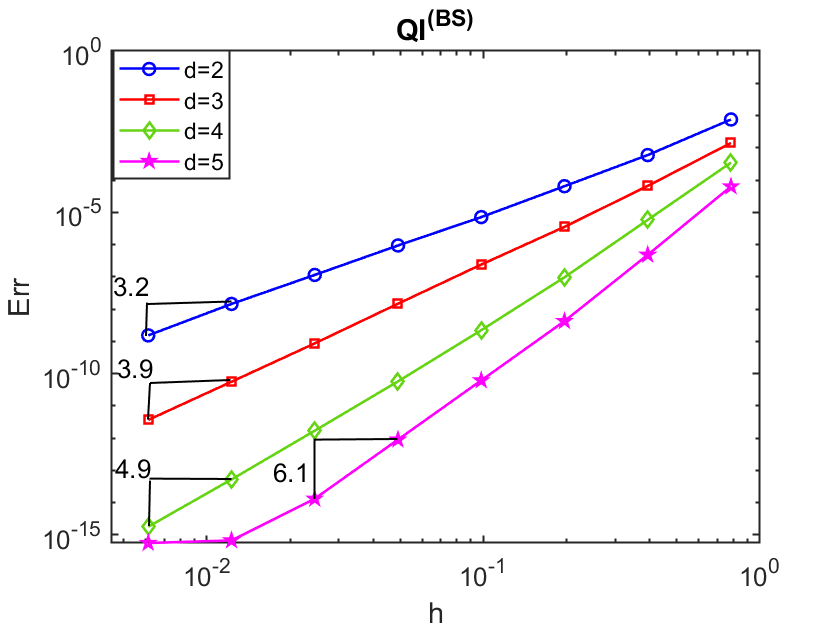}
    \includegraphics[height=5.15cm]{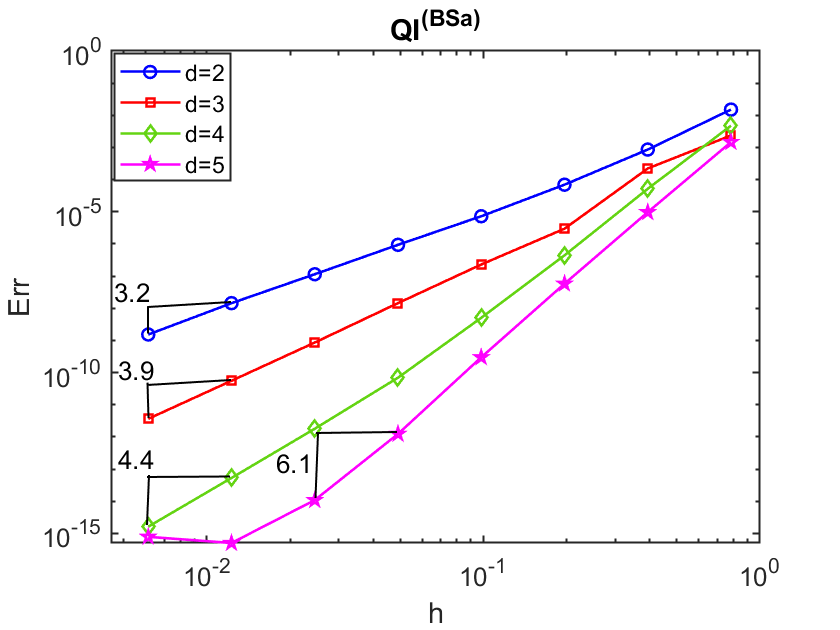}
    \caption{Convergence plots for the function $f(x)=\sin(x)$, varying the adopted QI degree.}
    \label{fig:QIex3}
\end{figure}
\begin{figure}
    \centering
    \includegraphics[height=5.15cm]{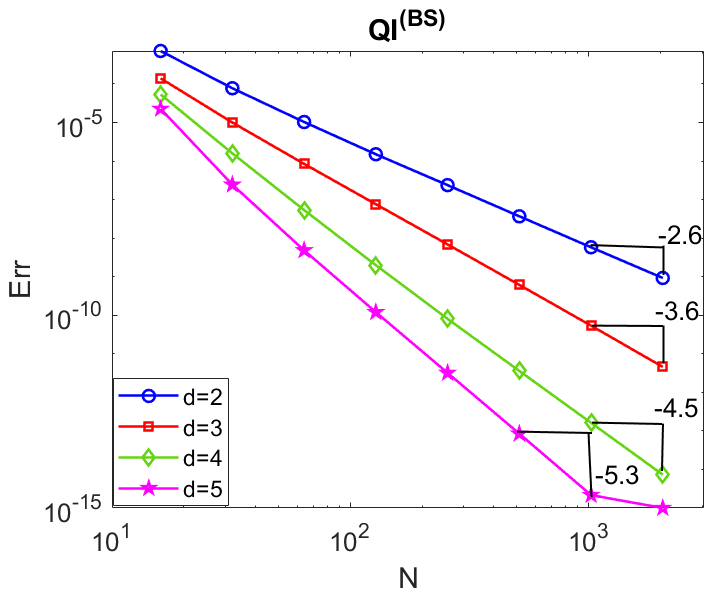}
    \includegraphics[height=5.15cm]{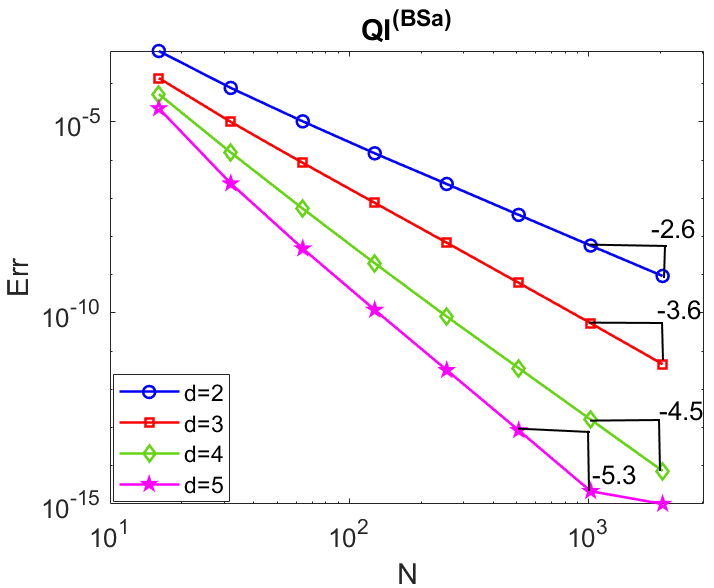}
    \caption{Convergence plots for the function $f_2$ with a non uniform knot partition, varying the used QI degree.}
    \label{fig:Qinonuni}
\end{figure}

\section{Numerical Tests}\label{test}

\subsection{Functions of One Variable}\label{test_una_dimensione}

We report the convergence results for the 
BSH Quasi Interpolant on the test function
\begin{equation}
  f_1(x)=e^{-x}\sin(5\pi x)\quad  x\in\left[-1,1\right]
\end{equation}
from~\cite{Mazzia1}.
Numerical examples are performed using the~\MATLAB{} Toolbox~\QIBSHPP{}. 
 
The Table \ref{tab:test1dabis2} shows the errors and the estimated 
orders of convergence of the Quasi Interpolant.  
For spline degree $d=3$ we apply to the function test the operators $Q_d^{(BS)}$ and $Q_d^{(BSa)}$, defined respectively in \eqref{eq:spline1D} and \eqref{eq:QIBSHa1D}. 
For the sake of comparison, we use also a discrete 
Quasi Interpolant of the same degree $d$, here denoted as $QI S$,
described in~\cite{SAB2} and the spline interpolant~\emph{spapi} 
from \MATLAB{}~\cite{MATLAB:2012}.
As usual, $N$ is the number of mesh steps between uniformly 
spaced $N+1$ spline knots, and errors are estimated by 
the infinity norm on  $1000$ points uniformly spaced 
in the domain of the functions, the timing  $\mathbf{t}$ 
is computed as the averaged time over $40$ runs (in secs.) 
for approximating the test functions and evaluating the 
spline using an Intel Core   I7-6500U $2.50$GHz with \MATLAB{} R2020a for Windows ($64$ bit).
The error values confirm the expected rate of convergence 
for $Q_d^{(BS)}$ and $Q_d^{(BSa)}$, as $d+1$.
 We also report the convergence behaviour
for the operator $QI^{(BS)}$ and $QI^{(BSa)}$ for $d$ ranging between $2$ and $5$ in Figure \ref{fig:QIex3}. 
Moreover, since the described quasi-interpolant operator 
can also be applied with a non-uniform knot partition,
in Figure \ref{fig:Qinonuni} we show the convergence behavior
for this case with the function $f_2 = \dfrac{e^{-x/s}-e^{(x-2)/s}}{1-e^{-2/s}  }$, and $s = 10^{-3/2}$.


Note that as $N$ increases, the error using $Q_d^{(BS)}$ and $Q_d^{(BSa)}$ 
becomes similar, especially for the choice $l=d+1$. 
To preserve the order of convergence it is enough to choose $l=d$, 
but the  part of the error due to the derivative approximation 
is dominant in this case.
These results suggest that the use of $Q_d^{(BSa)}$ can be 
competitive with other DQI interpolants, especially when
the approximation of higher order derivatives is required.

\begin{table}[tbp]
	\tbl{Convergence Analysis : function $f_1$, spline degree $3$, GBDF order
		$4$.\label{tab:test1dabis2}} {
		\centering
		\begin{tabular}{|c|c|c|c|c|c|c|c|c|c|c|c|c|}
			\hline
			\multicolumn{1}{|c|}{} &
			\multicolumn{3}{c|}{\QIBSH{}}&
			\multicolumn{3}{c|}{\QIBSH{}a} \\ \hline
			$\mathbf{N}$ &
			$\mathbf{Error}$ &
			$\mathbf{Order}$&
			$\mathbf{t}$&
			$\mathbf{Error}$ &
			$\mathbf{Order}$&
			$\mathbf{t}$
			\\ \hline
			\hline
			16 &
			\numprint{1.7e-1} & $**$ & \numprint{1.4e-5} &
			\numprint{8.5e-1} & $**$ & \numprint{2.4e-5}
			\\ \hline
			32 & 
			\numprint{8.8e-3} & 4.3 & \numprint{1.6e-5} &
			\numprint{7.3e-3} & 6.9 & \numprint{2.7e-5}
			\\ \hline
			64 & 
			\numprint{3.1e-4} & 4.8 & \numprint{1.4e-5} &
			\numprint{1.5e-4} & 5.6 & \numprint{2.0e-5}
			\\ \hline
			128 &
			\numprint{2.4e-5} & 3.7 & \numprint{1.4e-5} &
			\numprint{1.8e-5} & 3.1 & \numprint{2.4e-5}
			\\ \hline
			256 &
			\numprint{9.1e-7} & 4.7 & \numprint{1.8e-5} &
			\numprint{8.3e-7} & 4.4 & \numprint{3.6e-5}
			\\ \hline
			512 &
			\numprint{7.5e-8} & 3.6 & \numprint{2.8e-5} &
			\numprint{7.4e-8} & 3.5 & \numprint{6.0e-5}
			\\ \hline
			1024 &
			\numprint{4.2e-9} & 4.2 & \numprint{5.3e-5} & 
			\numprint{4.2e-9} & 4.1 & \numprint{1.1e-4}
			\\ \hline
			\hline
			\hline
			\multicolumn{1}{|c|}{} &
			\multicolumn{3}{c|}{\MATLAB{}}&
			\multicolumn{3}{c|}{\bf QIS}
			\\ \hline
			$\mathbf{N}$ &
			$\mathbf{Error}$ &
			$\mathbf{Order}$ &
			$\mathbf{t}$&
			$\mathbf{Error}$ &
			$\mathbf{Order}$&
			$\mathbf{t}$
			\\
			\hline
			\hline
			16 &
			\numprint{5.2e-1} & $**$ & \numprint{6.5e-4} &
			\numprint{4.9e-1} & $**$ & \numprint{1.3e-4}
			\\ \hline
			32 & 
			\numprint{1.1e-2} & 5.6 & \numprint{1.2e-3} & 
			\numprint{5.1e-2} & 3.3 & \numprint{1.3e-4}
			\\ \hline
			64 & 
			\numprint{2.6e-4} & 5.3 & \numprint{2.3e-3} &
			\numprint{4.0e-3} & 3.7 & \numprint{1.9e-4}
			\\ \hline
			128 &
			\numprint{2.2e-5} & 3.6 & \numprint{3.8e-3} &
			\numprint{2.7e-4} & 3.9 & \numprint{1.3e-4}
			\\ \hline
			256 &
			\numprint{9.0e-7} & 4.6 & \numprint{6.1e-3} &
			\numprint{1.6e-5} & 4.1 & \numprint{1.5e-4}
			\\ \hline
			512 &
			\numprint{7.5e-8} & 3.6 & \numprint{1.4e-2} &
			\numprint{1.1e-6} & 3.9 & \numprint{1.5e-4}
			\\ \hline
			1024 &
			\numprint{4.2e-9} & 4.1 & \numprint{3.3e-2} &
			\numprint{6.5e-8} & 4.0 & \numprint{1.8e-4}
			\\ \hline
			\hline
		\end{tabular}
	}
\end{table} 

Moreover, comparing the computational time (expressed in secs.) 
in the Table~\ref{tab:test1dabis2} and in Figure~\ref{fig:tempi1D}
top for $d=3$ and in Figure~\ref{fig:tempi1D} bottom,
for $d=5$, confirms the efficiency of the QI BSH
with respect to the other methods.

\subsection{Numerical solution of Boundary Value Problems}\label{sec:BVP}
The main aim of the library~\QIBSH{} is to improve the performance
of the~\MATLAB{} code~\PACKAGE{TOM} for the numerical solution 
of Boundary Value Problems (BVPs) that are assumed to have the first order system form,
\begin{equation}\label{bvp}
  y^\prime(x) = f (x,y), \quad a \leq x \leq b,
\end{equation}
where $y \in R^{m}, f : R \times  R^{m} \rightarrow R^{m}$,
with boundary conditions,
\[ g(y(a),y(b)) = 0. \]

\begin{figure}[!htb]
    \centering
    \includegraphics[width=6.5cm]{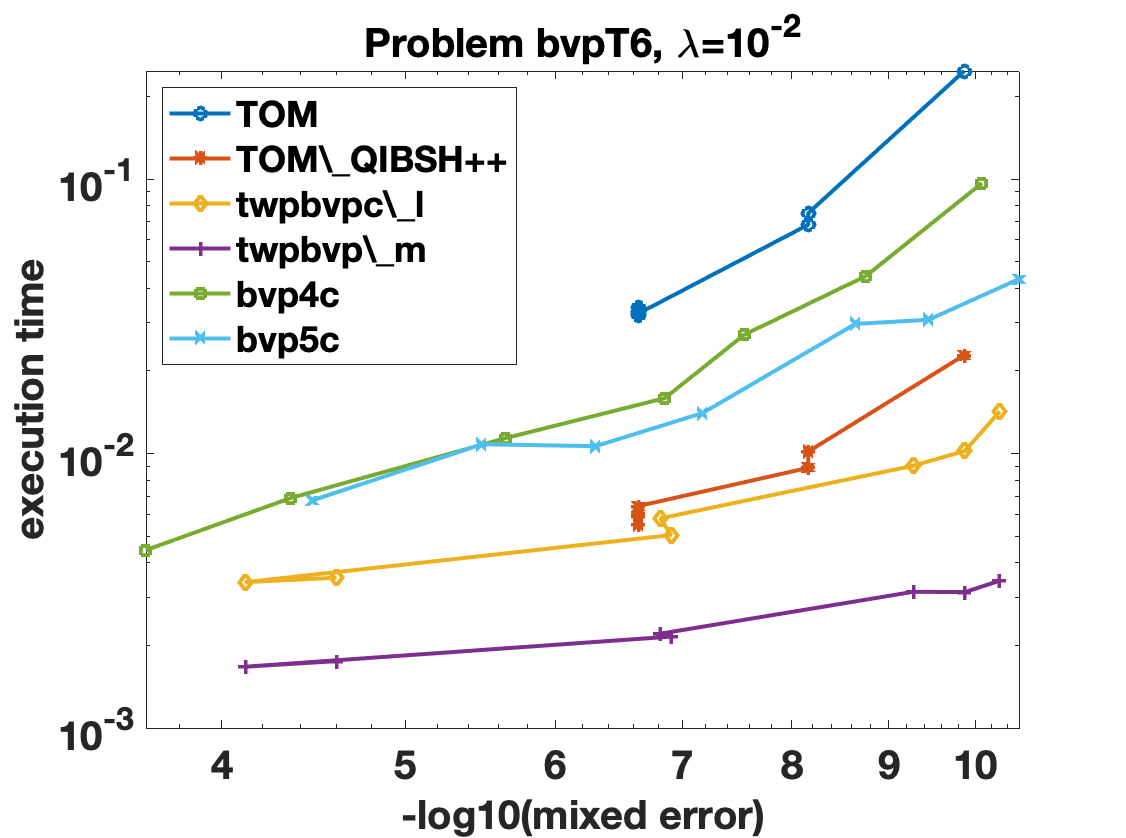}
    \includegraphics[width=6.5cm]{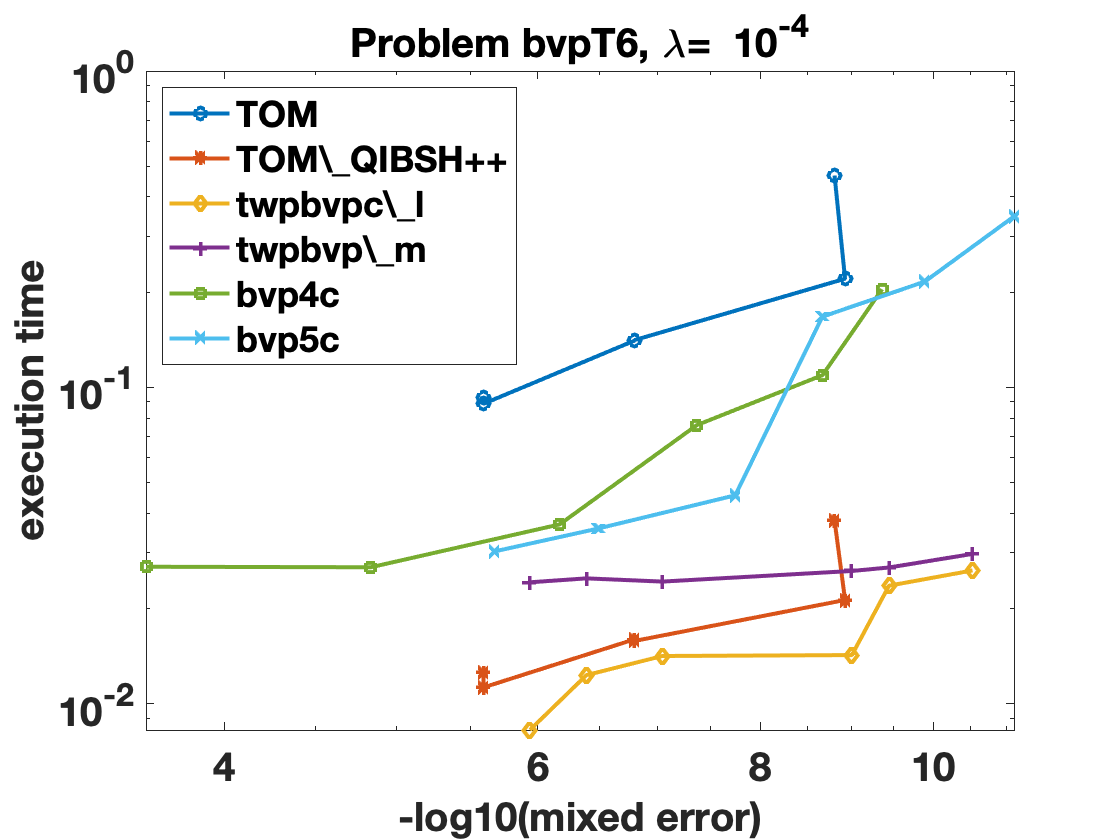}
    \caption{Work-precision diagram. Linear problem \PACKAGE{bvpT6} with different values of $\lambda$.}
    \label{fig:bvpT6}
    \centering
    \includegraphics[width=6.5cm]{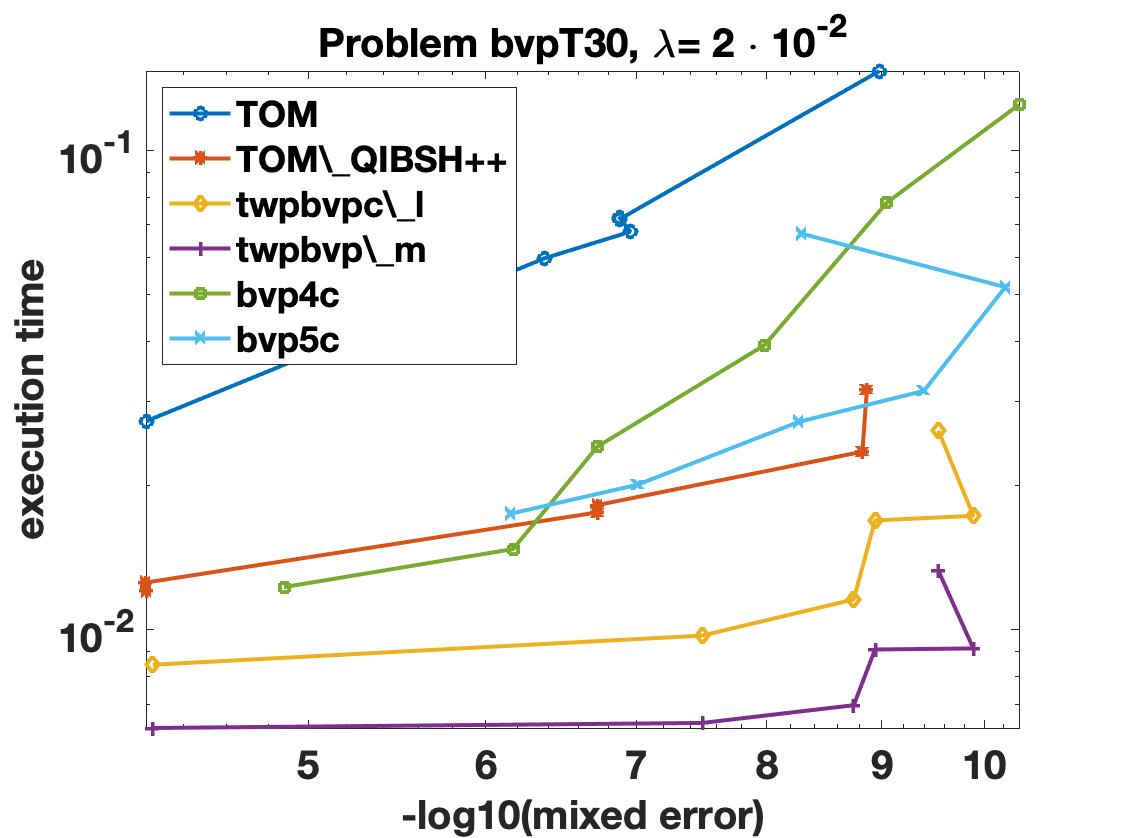}
    \includegraphics[width=6.5cm]{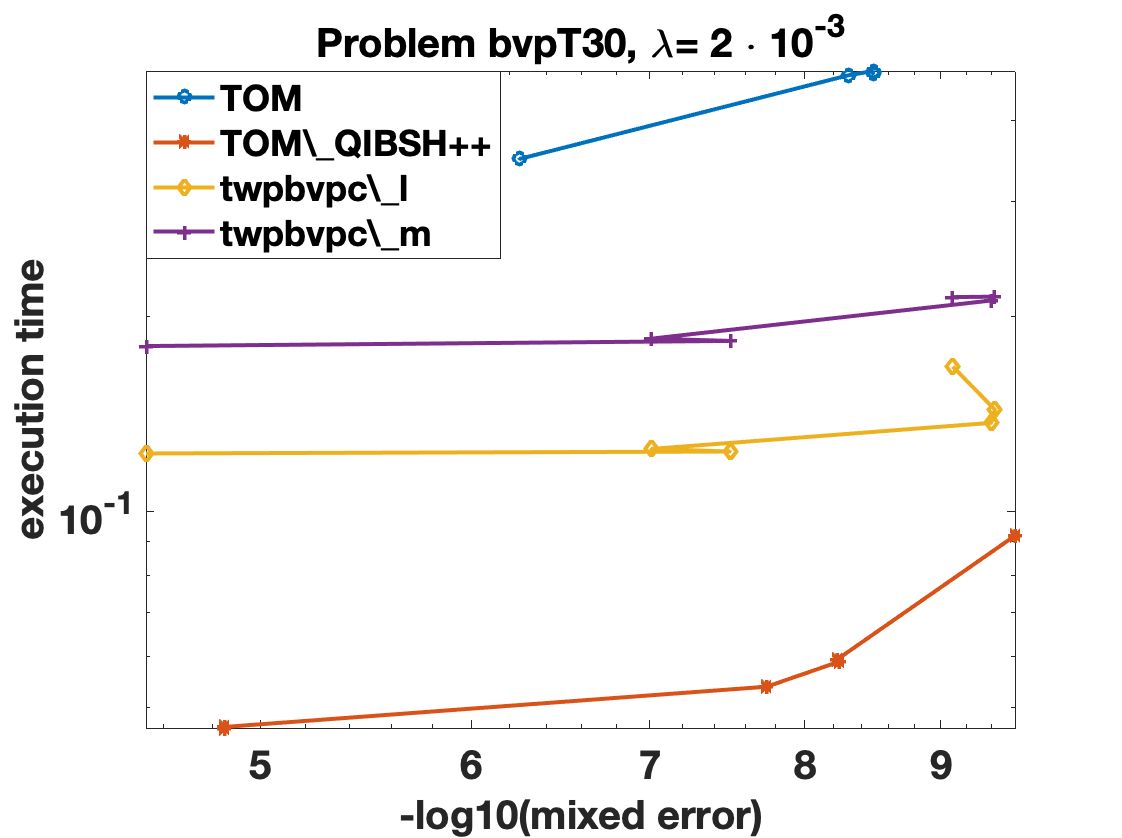}
    \caption{Work-precision diagram. Non linear problem \PACKAGE{bvpT30} with different values of $\lambda$.}
    \label{fig:bvpT30}
\end{figure}
The first release of the code has been described in~\cite{Mazzia.Trigiante.Mesh.04} whereas 
the update release that include the BS  linear multistep method 
is described in~\cite{Mazzia2009}. One of the main characteristic of the code is that it implements an hybrid mesh selection based on conditioning.
With the inclusion of the BS-scheme the code was able to solve
very difficult singularly perturbed BVPs, giving in output a 
continuous approximation of the solution, but  was not able to
have an efficient execution time for general problem, the Matlab
code spending most of the  time in computing the variable 
coefficients of the linear multistep methods.
We do not describe here the new release of the code but we 
just show how the use of the library~\QIBSH{} make this code faster with respect to the original~\MATLAB{} version and
competitive with the~\MATLAB{} codes \PACKAGE{bvp4c} and \PACKAGE{bvp5c}
and the code \PACKAGE{bvptwp}~\cite{cash2013algorithm}. 
The library has been integrated with the code TOM for both the
evaluation of the coefficients of the linear multistep BS 
method and the evaluation of the continuous extensions, 
needed for nonlinear problems when the mesh is changed.
The quasi-interpolant is also used  when the multistep
method is the Top Order Method, or another class of 
Boundary Value Method, as quasi-interpolation scheme.

Here we report some numerical experiment on two singularly
perturbed boundary value problems used in~\cite{cash2013algorithm}
and available in \cite{MazziaCashAIP2015}.
We choose two examples in this class that require
changes of the mesh  in order to compute the numerical solution and need a variable step-size. 
All the examples in this section have been run using~\MATLAB{} R2021b 
on iMac with a 3.6 GHz Intel core i9, 10 core.
 
The first numerical test chosen is the problem \PACKAGE{bvpT6},
the second is the nonlinear problem \PACKAGE{bvpT30}. Both are  singular perturbed problems  with a turning point. 
 
 We use different values of the parameter $\lambda$ and we select  input tolerances $atol=rtol$ as $10^{-3},10^{-4}, \cdots, 10^{-8}$. The code \PACKAGE{TOM} is used with the BS method of order 6 and the hybrid mesh selection denoted NSSE \cite{Mazzia.Ferrara.2022}, designed for the solution of non stiff problems.   The work precision diagrams  are reported in Figures \ref{fig:bvpT6}-\ref{fig:bvpT30}. The two codes called \PACKAGE{TOM}  and \PACKAGE{TOM\_QIBSH++}  denote the code \PACKAGE{TOM}  with and without the use of the library~\QIBSH{}, we compare them with the codes \PACKAGE{bvp4c}, \PACKAGE{bvp5c},
\PACKAGE{twpbvpc\_l}, \PACKAGE{twpbvpc\_m} and \PACKAGE{twpbvp\_m}. The last are some of the available codes in the package \PACKAGE{bvptwp}, the first one use the conditioning in the mesh selection and the Lobatto schemes, the second is based on Monoimplicit Runge-Kutta schemes. It is interesting to see that in all the experiments the use of the  library~\QIBSH{} reduce considerable the time and make the code \PACKAGE{TOM} comparable with the other available codes. We observe that for the nonlinear problem bvpT30 with parameter $\lambda = 10^{-3}$ only the codes with a mesh selection based on conditioning are able to compute a solution. In this case Figure \ref{fig:bvpT30} clearly show the efficiency of the code \PACKAGE{TOM\_QIBSH++}. Further experiments using the code can been found in \cite{Mazzia.Settanni.2021,Mazzia.Ferrara.2022}.

\subsection{Functions of Two Variables}\label{sec:2dimensioni}


We use the BSH tensor product on the well known~\emph{Franke} 
function from the test suite in \cite{FR1,RENKA792}.  

\begin{table}[htbp]
\tbl{Convergence Analysis for test function Franke:  bivariate tensor product spline of degrees $3\times3$ GBDF order $4$.\label{tab:test2dfranke2}}{
\centering
\begin{tabular}{|c|c|c|c|c|c|c|}
\hline
\multicolumn{1}{|c|}{}&
\multicolumn{3}{c|}{\QIBSH{}}&
\multicolumn{3}{c|}{\QIBSH{}a} \\ \hline
$\mathbf{N}$ &
$\mathbf{Err}$ &
$\mathbf{Ord}$&
$\mathbf{t}$&
$\mathbf{Err}$ &
$\mathbf{Ord}$&
$\mathbf{t}$
\\ \hline
\hline
  16 & 
  \numprint{2.9e-3} & $***$ & \numprint{9.6e-4} &
  \numprint{1.8e-3} & $***$ & \numprint{9.7e-4}
  \\ \hline
  32 &
  \numprint{1.1e-4} & 4.7 & \numprint{9.8e-4} &
  \numprint{1.3e-4} & 3.9 & \numprint{1.0e-3}
  \\ \hline
  64 &
  \numprint{5.2e-6} & 4.4 & \numprint{1.1e-3} &
  \numprint{2.2e-6} & 5.9 & \numprint{1.1e-3}
  \\ \hline
 128 &
 \numprint{2.9e-7} & 4.1 & \numprint{1.4e-3} &
 \numprint{2.4e-7} & 3.2 & \numprint{1.5e-3}
 \\ \hline
 256 & 
 \numprint{1.6e-8} & 4.2 & \numprint{2.8e-3} &
 \numprint{1.5e-8} & 4.0 & \numprint{3.6e-3}
 \\ \hline
 512 &
 \numprint{1.1e-9} & 3.9 & \numprint{1.1e-2} &
 \numprint{1.1e-9} & 3.8 & \numprint{1.2e-2}
 \\ \hline
1024 &
\numprint{7.2e-11} & 3.9 & \numprint{4.5e-2} &
\numprint{7.2e-11} & 3.9 & \numprint{5.7e-2}
\\ \hline

\multicolumn{1}{|c|}{}&
\multicolumn{3}{c|}{\MATLAB{}}&
\multicolumn{3}{c|}{QIS}
\\ \hline
$\mathbf{N}$ &
$\mathbf{Err}$ &
$\mathbf{Ord}$&
$\mathbf{t}$&
$\mathbf{Err}$ &
$\mathbf{Ord}$&
$\mathbf{t}$
\\ \hline
\hline
  16 & 
  \numprint{2.1e-3} & $***$ & \numprint{4.4e-3} &
  \numprint{9.3e-3} & $***$ & \numprint{2.3e-3}
  \\ \hline
  32 &
  \numprint{7.9e-5} & 4.7 & \numprint{6.6e-3} &
  \numprint{7.9e-4} & 3.6 & \numprint{2.3e-3}
  \\ \hline
  64 &
  \numprint{4.6e-6} & 4.1 & \numprint{9.5e-3} &
  \numprint{5.2e-5} & 3.9 & \numprint{3.2e-3}
  \\ \hline
 128 &
  \numprint{2.9e-7} & 4.0 & \numprint{1.8e-2} &
  \numprint{3.4e-6} & 3.9 & \numprint{3.4e-3}
  \\ \hline
 256 &
  \numprint{1.6e-8} & 4.2 & \numprint{3.2e-2} &
  \numprint{2.2e-7} & 4.0 & \numprint{7.8e-3}
  \\ \hline
 512 &
  \numprint{1.1e-9} & 3.9 & \numprint{9.4e-2} &
  \numprint{1.4e-8} & 4.0 & \numprint{2.6e-2}
  \\ \hline
1024 &
  \numprint{4.1e-14} & 14.7 & \numprint{3.7e-1} &
  \numprint{8.5e-10} & 4.0 & \numprint{8.5e-2}
  \\ \hline

\end{tabular}
}
\end{table} 
The Table \ref{tab:test2dfranke2} summarizes the results on the test function for the tensor product spline of degrees $d_x=d_y=3$, and the choice $l=d_x+1=d_y+1$ for the  order of the finite difference scheme. Again, we are 
comparing the behavior of the tensor product formulation of the operators 
$Q_d^{(BS)}$ and $Q_d^{(BSa)}$, to the tensor product of the DQI $QIS$, and to the one from \MATLAB{} \emph{spapi}.  
Infinity norm errors are computed against the exact values on a uniform grid of $101 \times 101$ points in the unit square.
The numerical experiments are carried out on a personal computer Intel Core   I7-6500U 2.50GHz with \MATLAB{} 2020a for Windows (64 bit).

\begin{figure}[!htb]
    \centering
    \includegraphics[width=6.5cm]{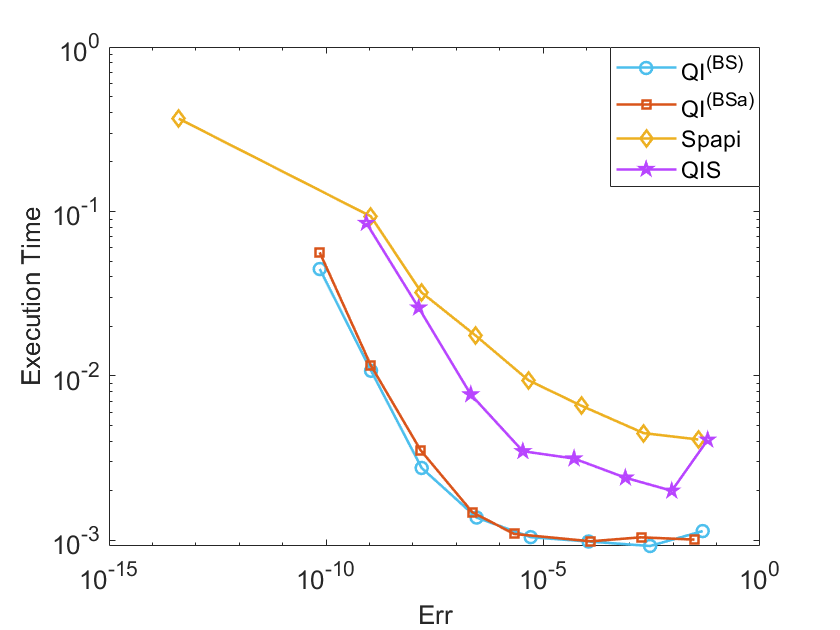}
    \includegraphics[width=6.5cm]{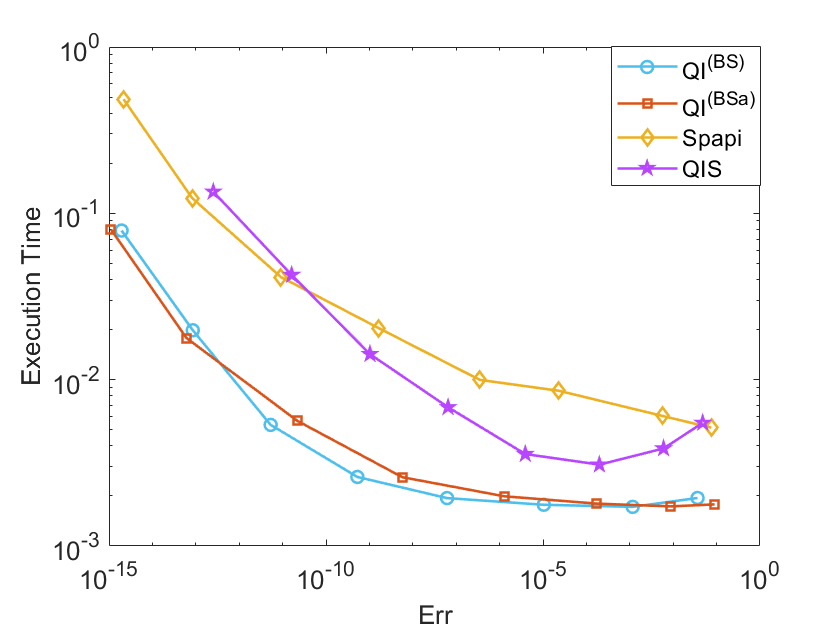}
    \caption{Work precision diagram for the Franke function. On the left the used bidegree $d=3$. On the right, the used bidegree $d=5$.}
    \label{fig:TimeFranke}
\end{figure}

When fewer knots are available, it is convenient to use the operator
$Q_d^{(BSa)}$, choosing at least $l=d$\footnote{ since in the following tests we always adopt $d_x$=$d_y$, we simplify the notation by using the letter $d=d_x=d_y$ and  denoting  the used tensor product of  bidegree $d$.}
otherwise the error for the first derivatives approximations gets larger. On the 
contrary, as $N$ increases, and a larger number of knots is available, the user 
can use a finite difference scheme increasing the approximation order $l$. Indeed, for $N\geq 256$,  the approximation for BSH  gives the same results. 
In particular they tend to have the same behavior of BSH, as if the partial 
derivatives  values would have been available. The Hermite Quasi Interpolant in 
these cases has smoothed all the first partial derivative approximating errors. 
Note that the computational time between $Q^{(BS)}$ and $Q^{(BSa)}$ does not 
differ so much. In fact, this is due to the 
efficient approximation for the derivatives, which is not adding any significant 
computational cost to the one of the QI.\\
Indeed, comparing the time efficiency of the~\QIBSHPP{} tensor product
operator to  the other interpolants tested, confirms its good behaviour as an approximating method for functions of two variables in terms 
of goodness of fit and timing.  In Figure \ref{fig:TimeFranke} we also report the work precision diagrams for the Franke function with biddegree $d=3$ and $d=5$.
This motivates once more  the use of the Hermite Quasi Interpolant
even when partial derivatives are not directly available,
a situation often occurring in real applications, as it will be shown later.

\begin{figure}[!htb]
    \centering
    \includegraphics[width=6.5cm]{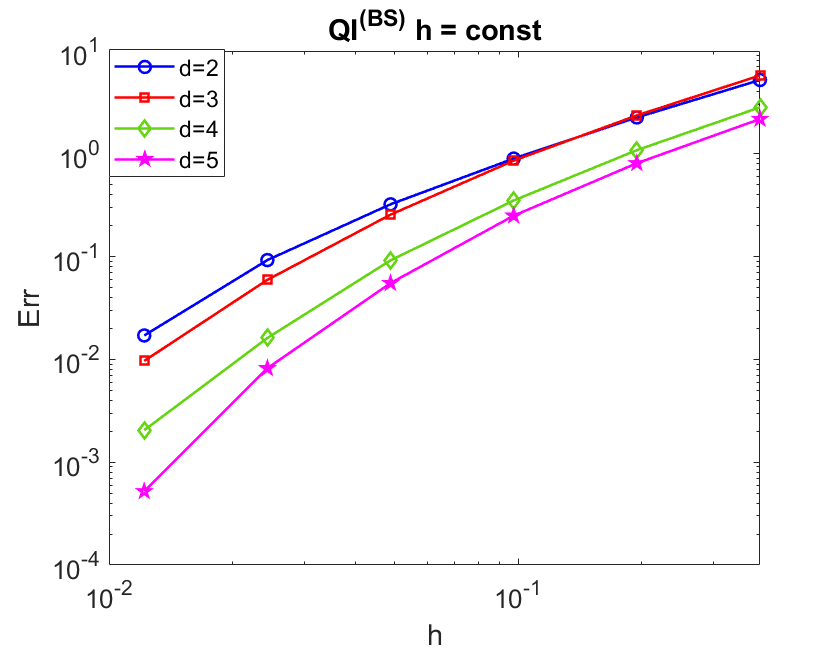}
    \includegraphics[width=6.5cm]{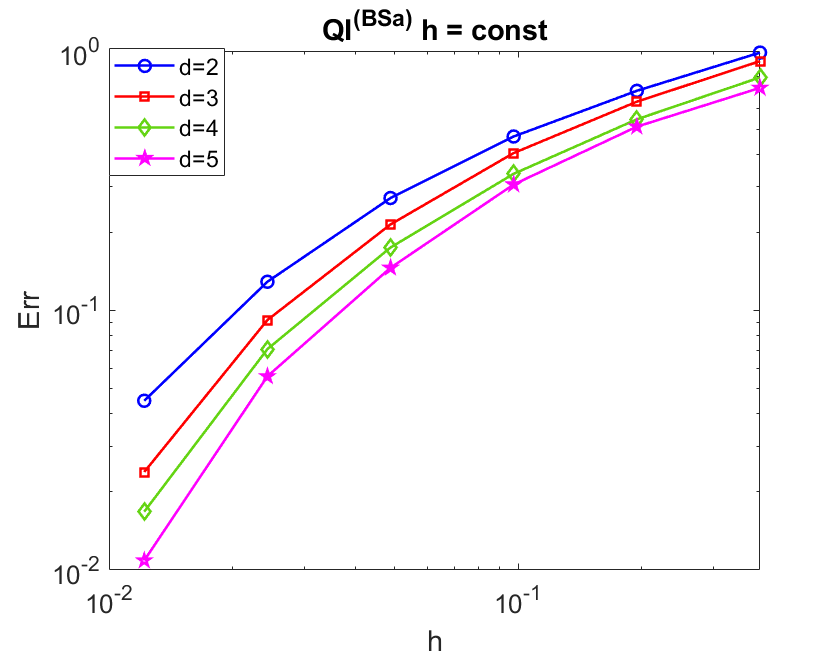}
    
    \includegraphics[width=6.5cm]{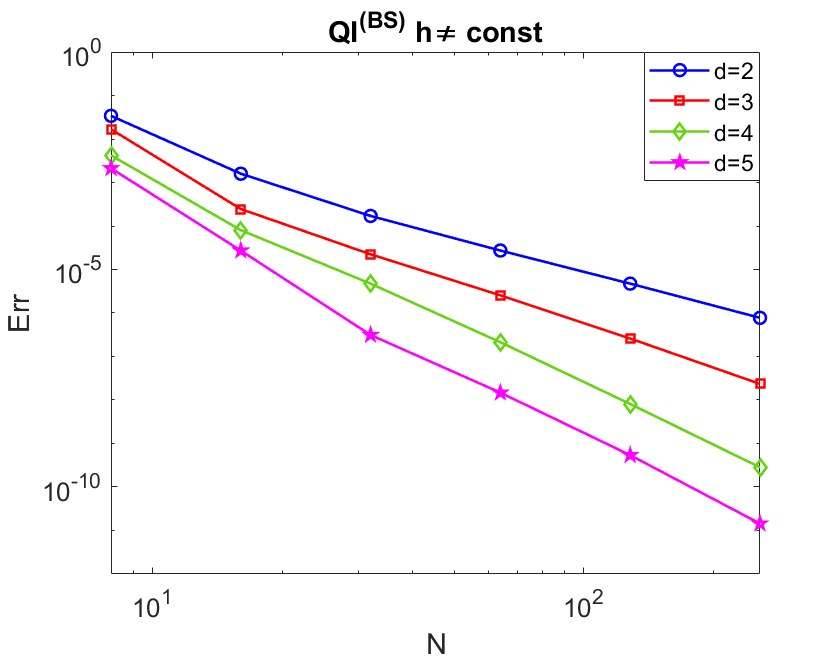}
    \includegraphics[width=6.5cm]{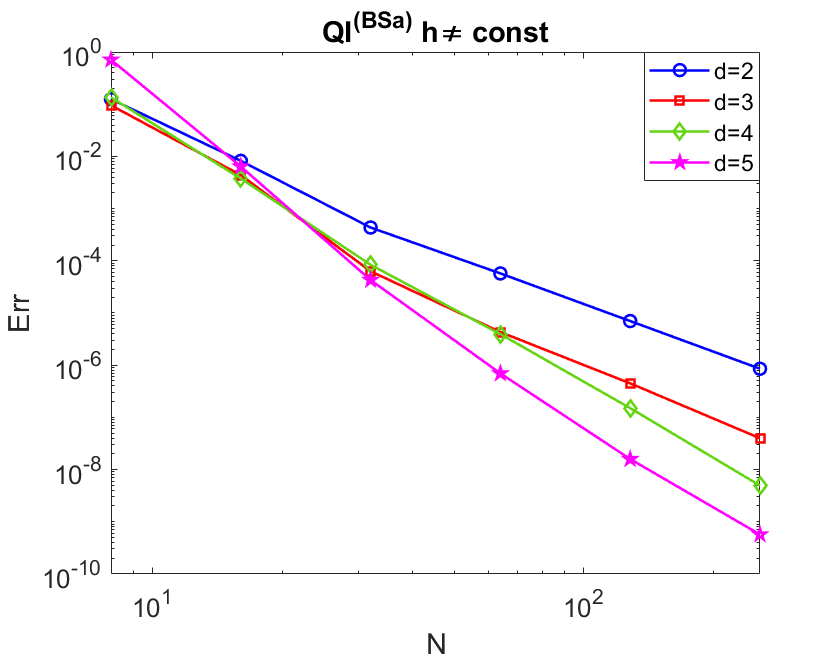}
    \caption{Convergence analysis for the approximation of the Scherk's first surface. On the top, the results obtained with the operators $QI^{(BS)}$ and $QI^{(BSa)}$ and a uniform knot partition. At the bottom, the achieved results with the same operators but a non uniform knot partition. In both cases the spline bidegree $d$ varies from $2$ to $5$.}
    \label{fig:ConvSherk}
\end{figure}

Also for the 2D case, we analyze an example where the convergence can benefit from a non uniform mesh partition. In particular, we approximate the Schrek's first surface $f(x,y) =\displaystyle{\log\left(\frac{\cos(y)}{\cos(x)}\right)}$ restricting our-self to the square domain $\displaystyle{\left[-\frac{\pi}{2}+\eta, \frac{\pi}{2}-\eta\right]^2}$, with $\eta = 10^{-2}$. 
\begin{figure}[!htb]
    \centering
    \includegraphics[width=6.5cm]{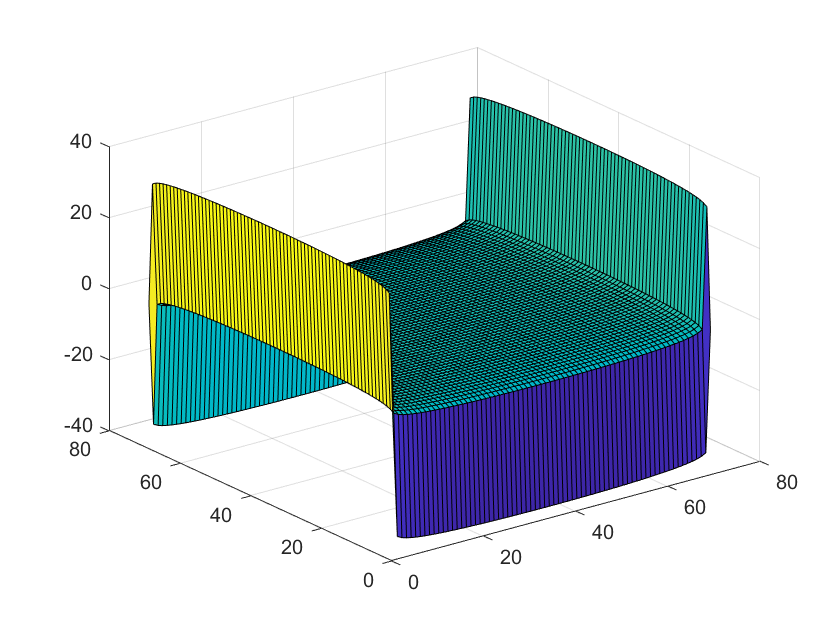}
    \caption{Scherk's first surface. }
    \label{fig:Sherk}
\end{figure}
As shown in Figure  \ref{fig:Sherk}, the considered example is a minimal area surface  which exhibits a high variation along each side of the definition domain. Therefore, a uniform knot partition will not guarantee a suitable error reduction at the boundary layers. The convergence results are shown in Figure \ref{fig:ConvSherk} also considering the spline bidegree $d$ varying from $2$ to $5$.

\subsection{Surface parameterization with high smoothness for complex geometries}\label{sub:surface}
\begin{figure}[!htb]
    \centering
    \includegraphics[height=5cm, trim = 2.5cm 0 2cm 0, clip=false]{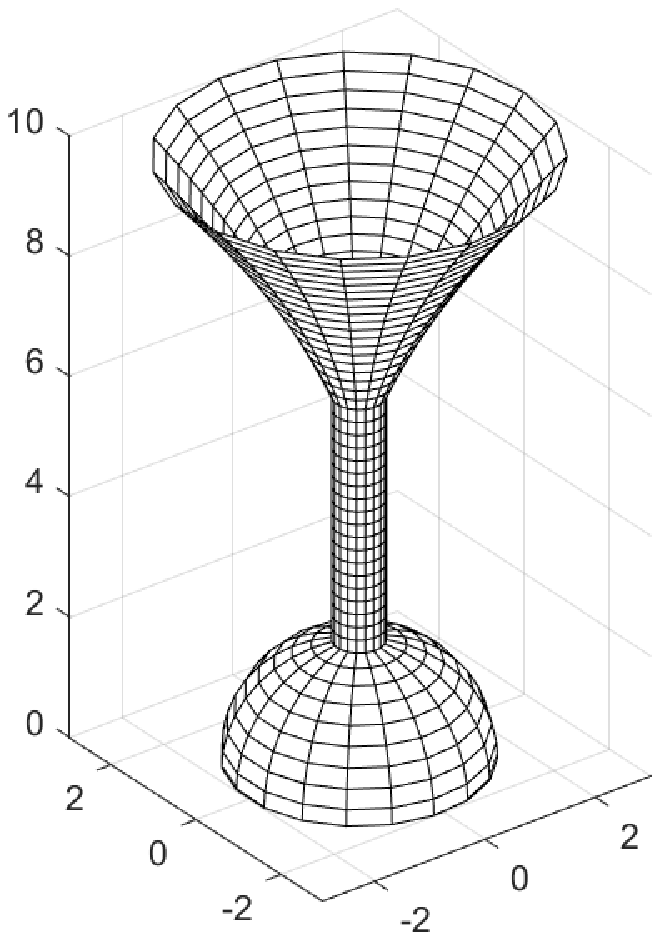}
    \includegraphics[height=5cm, trim = 2.5cm 0 2cm 0, clip]{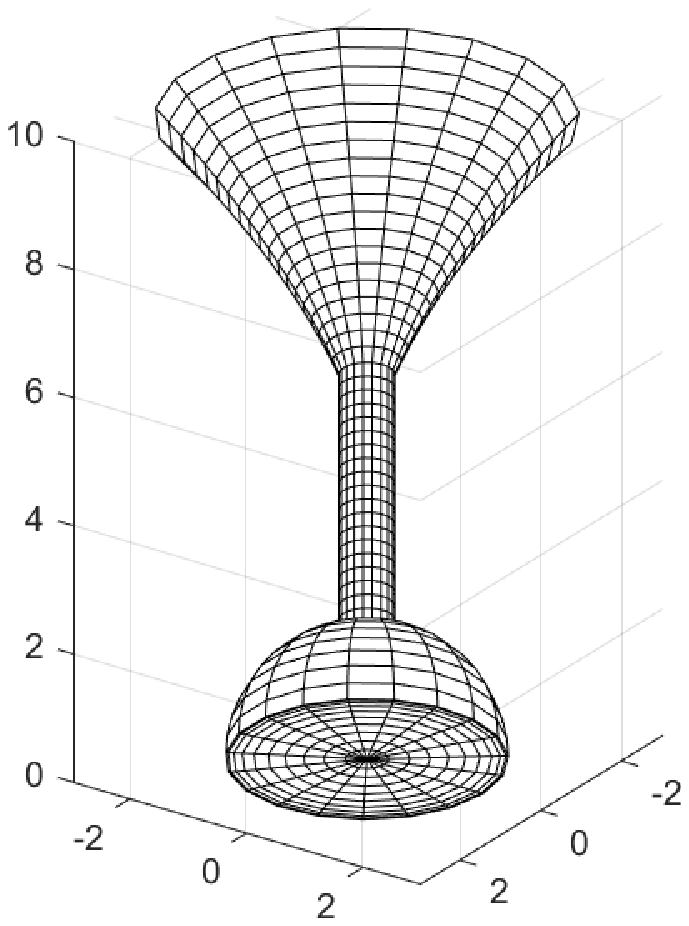}
    \caption{Initial surface. The primitive shapes of a hollow hemisphere, a cylinder and a conical frustum are joined with $C^0$ continuity.}
    \label{fig:multipatch}
\end{figure}

\begin{figure}[!htb]
    \centering
    \includegraphics[height=5cm]{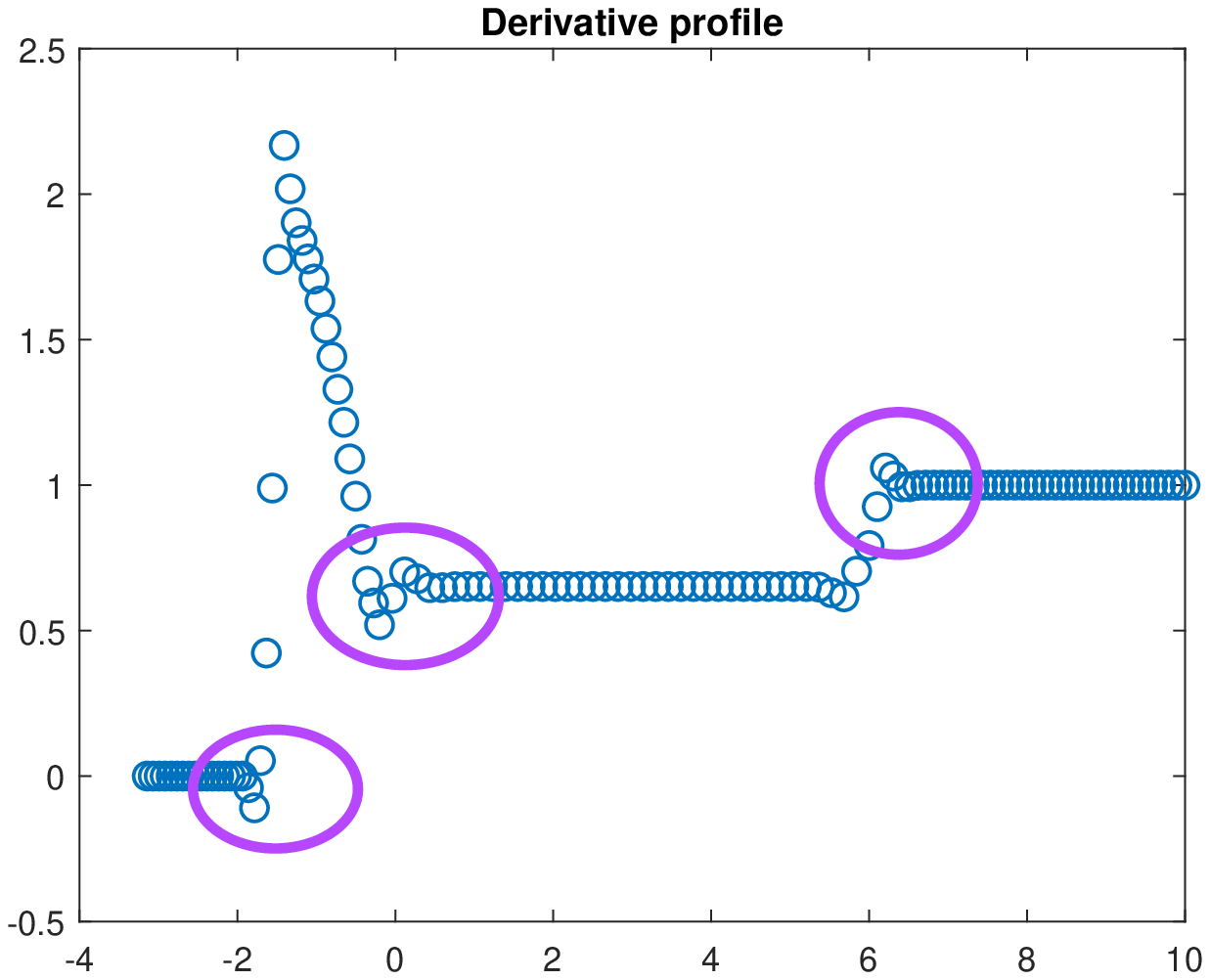}
    \includegraphics[height=5cm]{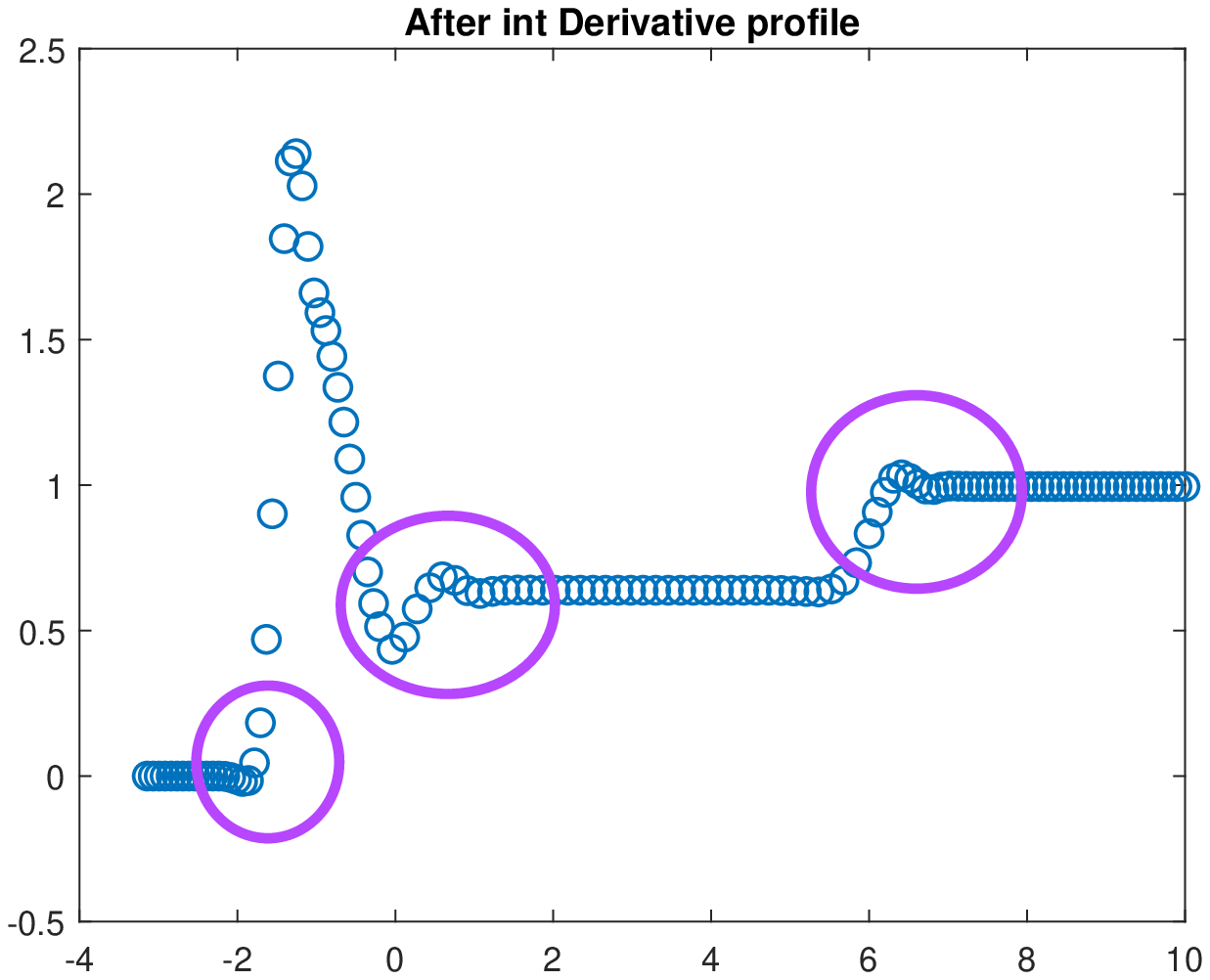}
    \caption{First derivative profile with respect to $y$ direction. Left: the original construction. Right: the regularized function after integration.}
    \label{fig:der-multi}
\end{figure}

In this example we demonstrate how the tensor product BSH QI can be useful to 
produce spline parameterizations of complex geometries with the desired smoothness in any direction. We construct a complex geometric model by assembling together 
three primitive shapes: hollow hemisphere, cylinder and conical frustum. The 
resulting object is a glass geometry shape, see Figure \ref{fig:multipatch}. In 
order to construct a continuous spline approximation of the considered object we 
need to reparametrize every shape in such a way to obtain conforming parameter 
domains. More in detail, we consider a hollow hemisphere $S$ with radius $r=2$ 
described by a parameterization $F_S(\theta,\phi): \Omega_S \rightarrow S$ with 
$\Omega_S:= [ 0, 2\pi] \times [-\pi, -0.2)$ with the following coordinates 
representation:
\begin{align*}
  X_S = & r\cos(\theta)\sin(-\phi);\\
  Y_S = & r\sin(\theta)\sin(-\phi);\\
  Z_S = & r\cos(\phi).
\end{align*}
We consider only the values for $Z_S>0$. 
The cylinder $C$ considered here can be described by the following $F_C(\theta,h):\Omega_C \rightarrow C$, with $\Omega_C = [0, 2\pi] \times [-0.2,6]$. The physical coordinates can be expressed as:
\begin{align*}
X_C = &~ r_C\cos(\theta);\\  
Y_C = &~ r_C\sin(\theta);\\
Z_C = &~ \frac{(6-r\cos(-0.2))}{(6+0.2)}(h+0.2) + r\cos(-0.2);
\end{align*}
with $r=2$ and $r_C = r\sin(0.2)$.
Finally the conical frustum $V$ is defined with $F_{V} (\theta, v): \Omega_{V}\rightarrow V$, with $\Omega_{V} = [0, 2\pi] \times (6,10]$. The physical coordinates are computed as:
\[
X_{V} = a\frac{v\cos(\theta)}{10};\qquad
Y_{V} = a\frac{v\sin(\theta)}{10};\qquad
Z_{V} = v;
\]
with $a$ uniformly varying from $|r_c h_f|/6$ to $3$.
The three shapes are parameterized in such a way that the resulting
geometries can be physically joined with $C^0$ continuity and 
their parameter domains can be assembled to form a unique domain
$\Omega$ defined as: $\Omega = \Omega_S \cup \Omega_C \cup \Omega_V$
where we define our quasi-interpolant spline approximation $F$.
In this case $\Omega = [0,\pi]\times [-\pi, 10]$. 
If we are interested in constructing a $C^0$ representation $F$, 
at this stage we only need to call the constructor for the 2D
object~\QIBSH{} and we can compute a spline representation $F$ 
with a chosen bidegree $d$, periodic along with the first direction,
by approximating the $F_x$, $F_y$ and $F_{xy}$ with a finite 
difference scheme of order $\ell = \max\{d_x,d_y \} +1$.
Although the resulting surface has in principle $C^{d-1}$ smoothness,
its derivatives might present sharp variation and/or unwanted 
oscillations, see Figure~\ref{fig:der-multi}, left, where we plotted 
the first derivative profile with respect to the $y$ direction.
To get an improved parametric representation, we therefore proceed as follows.
Firstly the desired smoothness for the final approximation $F$ is fixed.
For practical purposes we limit our-self to the case of constructing 
a $C^1$ spline parameterization.
Hence, the chosen degree should be $d_x,d_y\geq 2$. 
Secondly, since for this example the derivatives with respect 
to $x$ direction are almost zero, we will tackle only 
the derivatives with respect to the $y$ direction.

The proposed algorithm can be summarized with the following steps,
\begin{itemize}
    \item[(1)] A discrete approximation $\tilde{F}_{yy}$ of the second derivatives in the $y$ direction is constructed, by using finite centered differences. 
    \item [(2)] A continuous model of $\tilde{F}_{yy}$ is provided by applying the $Q^{(BSa)}$ operator, thus obtaining $\tilde{F}_{yy}\approx QA_{yy}$.
    \item [(3)] The first derivative $QA_y$ in $y$ direction is now computed as,
    $$
    QA_y(\hat{x},y_z):= \int_{c}^{y_z} QA_{yy}(\hat{x}, z) \, dz,
    $$
    where $y_z, \, z = 0, \ldots, M$, knots of the mesh.
\end{itemize}
The output of this procedure can be seen in Figure~\ref{fig:der-multi}, right.
We can visually appreciate how the function results more regularized.
Some oscillations are still visible, in fact, the steps (1)-(3) can
be performed by starting from any order of derivation and then can be backward iterated till the computation of the approximant surface for the original one.

\subsection{Functions of three Variables}\label{test_tre_dimensioni}

We consider the following volume:
\begin{align*}
f(x,y,z)=    \frac{\sqrt{64-81((x -1/2)^2 + (y -1/2)^2 + (z-1/2)^2)}}{9}-\frac{1}{2}\quad x,y,z \in [0,1].
\end{align*}
In Table~\ref{tab:test3d} we report the maximum error and the estimated order 
of convergence for the operator $Q_d^{(BSa)}$, by varying the used degree
$d$ for each direction and the number of mesh steps $N$ between 
uniformly spaced $N + 1$ knots. In the following test $d_x=d_y=d_z=d$
and the approximation order $\ell$ for the derivatives is computed 
as $d+1$ for $d$ odd and as $d+2$ for $d$ even, in order to produce
a symmetric output in line with the original function $f$.
The order of convergence $d+1$ is reached in all the considered cases.

\begin{table}[htbp]
	\tbl{Convergence Analysis for the volume $f$ in Section \ref{test_tre_dimensioni}. \label{tab:test3d}} {
		\centering
		\begin{tabular}{|c|c|c|c|c|c|c|c|c|}
			\hline
			\multicolumn{1}{|c|}{}&\multicolumn{2}{c|}{\QIBSH{}a $d=2$}&\multicolumn{2}{c|}{\QIBSH{}a $d=3$}&\multicolumn{2}{c|}{\QIBSH{}a $d=4$}&\multicolumn{2}{c|}{\QIBSH{}a $d=5$}  \\ \hline
			$\mathbf{N}$ &
			$\mathbf{Err}$ &
			$\mathbf{Ord}$&
			$\mathbf{Err}$ &
			$\mathbf{Ord}$&
			$\mathbf{Err}$ &
			$\mathbf{Ord}$&
			$\mathbf{Err}$ &
			$\mathbf{Ord}$
			\\ \hline
			\hline
			16 & 
			\numprint{2.3e-3} & $**$ &
			\numprint{1.7e-3} & $**$ & 
			\numprint{9.5e-4} & $**$ & 
			\numprint{9.6e-4} & $**$ 
			\\ \hline
			32 &
			\numprint{3.7e-4} & 2.6 & 
			\numprint{2.7e-4} & 2.7 & 
			\numprint{1.1e-4} & 3.1 & 
			\numprint{1.2e-4} & 3.1 
			\\ \hline
			64 &
			\numprint{2.8e-5} & 3.8 & 
			\numprint{1.6e-5} & 4.1 & 
			\numprint{4.4e-6} & 4.7 & 
			\numprint{4.3e-6} & 4.7 
			\\ \hline
			128 &
			\numprint{1.9e-6} & 3.9 & 
			\numprint{5.8e-7} & 4.8 & 
			\numprint{8.8e-8} & 5.6 & 
			\numprint{6.7e-8} & 6.0 
			\\ \hline
			256 & 
			\numprint{2.8e-7} & 2.8 & 
			\numprint{1.4e-8} & 5.4 & 
			\numprint{3.1e-9} & 4.9 & 
			\numprint{9.2e-10} & 6.2 
			\\ 
			\hline
		\end{tabular}
	}
\end{table}

\subsection{Real World Data Tests}\label{test_dati_reali}

In this Section we report two real applications where the use of the 
library~\QIBSH{} gave interesting results and improvements with respect to 
other usual techniques.

\subsubsection{Continuous Digital Elevation Models}

We consider data-sets available to produce Digital Elevation 
Models (DEM). DEM  is a digital model or a 3D representation 
of a terrain's surface altitude with respect to the mean sea level.
Technically DEM contains only the elevation information 
without taking into account possible vegetation, buildings, 
or other types of objects. DEM are generated by using the 
elevation information from points spaced either at regular 
or irregular intervals.
In the first case, when the points are collected in regular grids,
we talk about raster DEM, in the latter case, the points are 
arranged in triangular irregular networks, hence, we refer to vector DEM.
In this example we use the \PACKAGE{NASADEM}\footnote{NASA JPL (2021).
\PACKAGE{NASADEM} Merged DEM Global 1 arc second V001.
Distributed by \PACKAGE{OpenTopography}.
\url{https://doi.org/10.5069/G93T9FD9 Accessed: 2021-08-22 dataset}}:
a modernization  of Global Digital Elevation Models.
The satellite data are preprocessed according to several 
optimization and interpolation algorithms and they are provided
in grid form. \PACKAGE{NASADEM} products are freely available 
through the Land Processes Distributed Active Archive Center 
(LP DAAC) at $1$ arc-second spacing ( $\sim 30$ meters).
Data were collected from February 11, 2000 to February 22, 2000. 
In order to produce a continuous model of the discrete dataset,
since the given samples are uniformly spaced (raster DEM),
it is reasonable to adopt a tensor product approach.
In addition, the produced continuous model should still
be able to capture the abrupt changes in the terrain shape,
so it is reasonable to require up to  $C^2$ smoothness. 

Among the many application of DEM a continuous model might 
be used for 3D rendering visualization purposes, 
hydrological and geomorphological investigations, 
rectification of satellite imagery, terrain correction and so on. 

\begin{figure}[!htb]
    \centering
    \includegraphics[width=6cm]{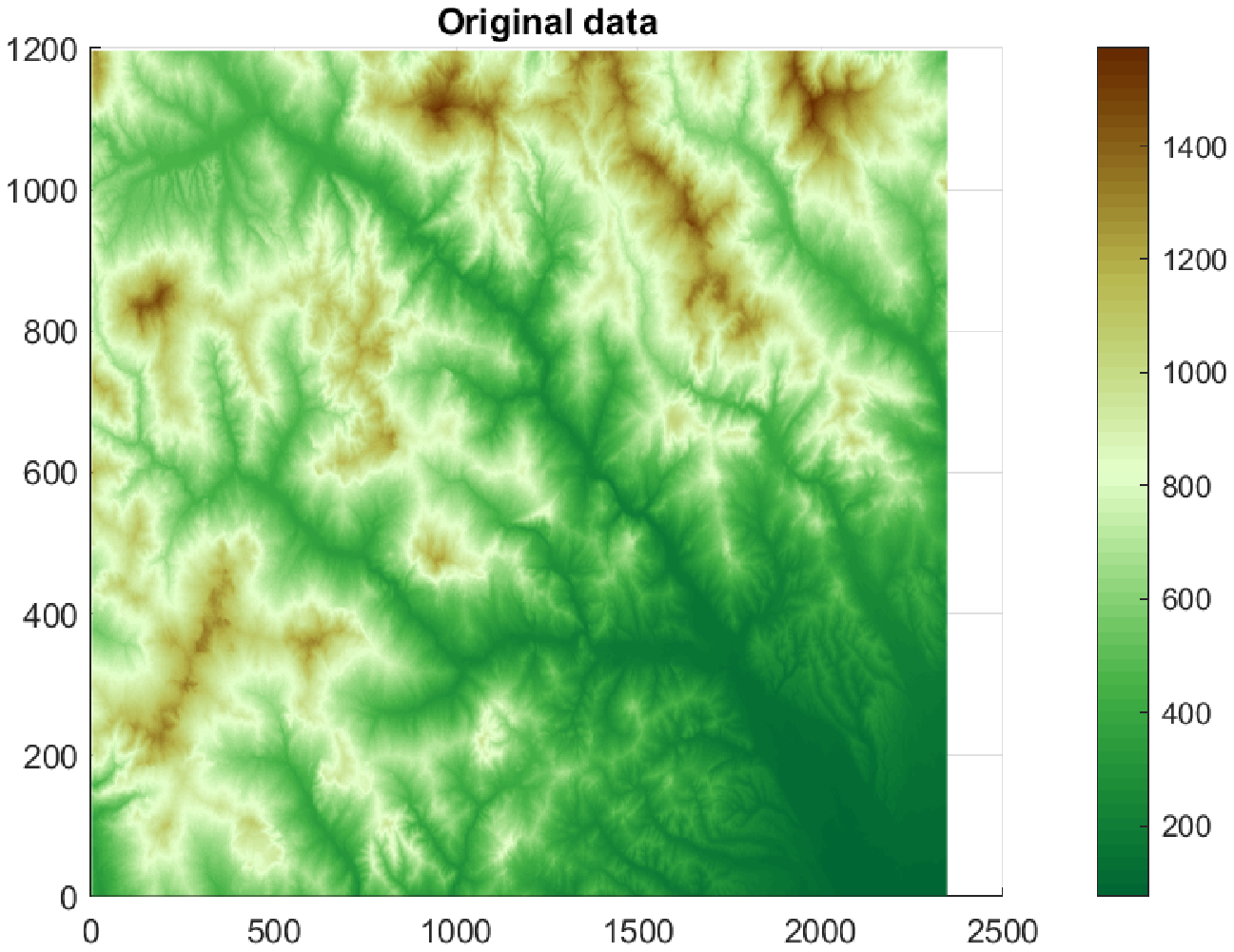}
    \includegraphics[width=6cm]{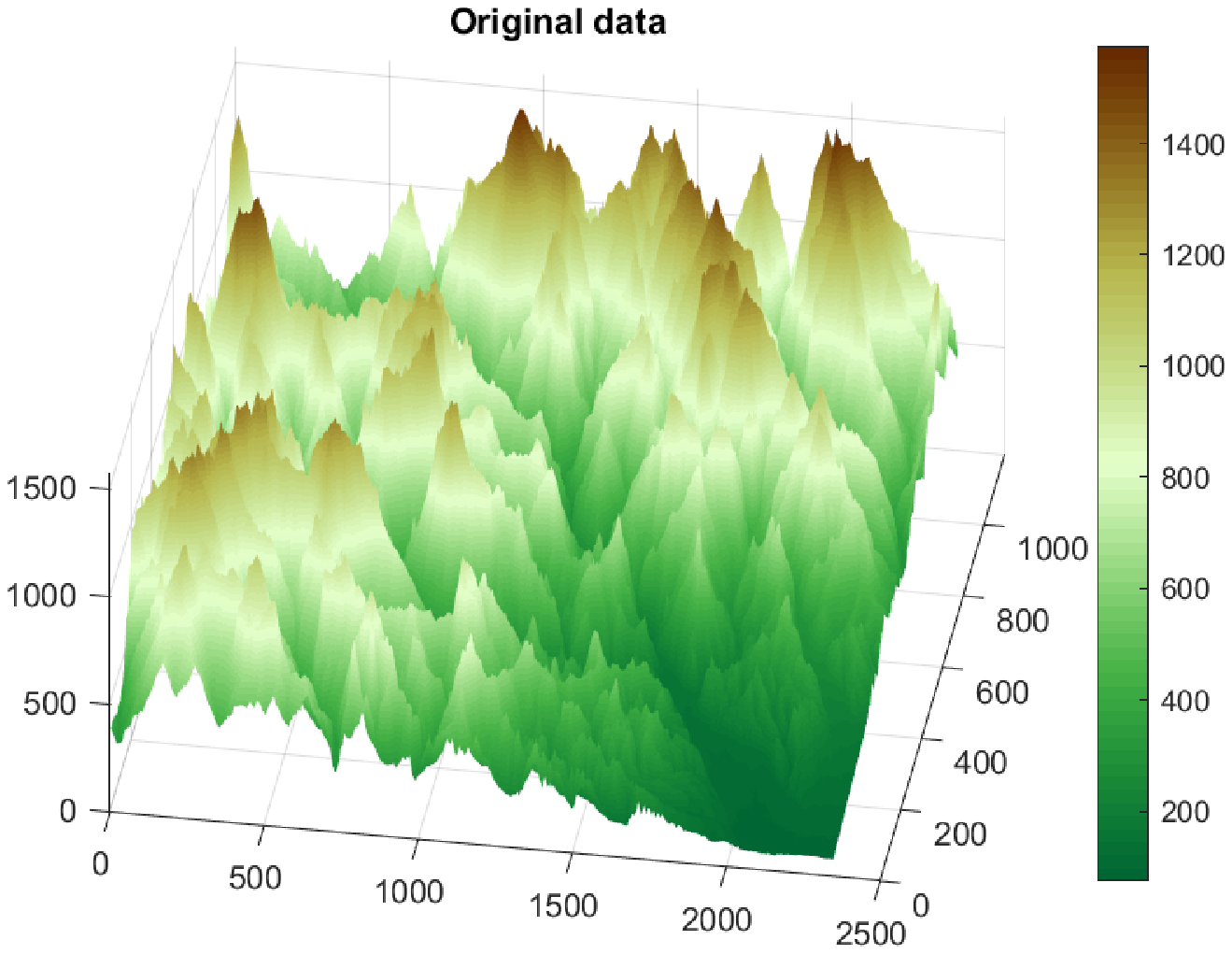}
    \caption{\PACKAGE{NASADEM} in the central Italian Appenini.}
    \label{fig:Nasadem}
\end{figure}

In this example we select a terrain matrix of size 
$1197\times 2347$ in the Italian Appenini mountain region, 
see Fig.~\ref{fig:Nasadem}. In Table~\ref{tab:QIterrain} we report
the root mean square error (RMSE) and the normalized RMSE (NRMSE)
on the quasi-interpolation nodes. 

\begin{table}[!htb]
\tbl{RMSE and NRMSE on the quasi-interpolation nodes.\label{tab:QIterrain}} {
\centering
\begin{tabular}{|c|c|c|}
\hline
 &{\QIBSH{}a}& {\QIBSH{}a} \\
 & $d=2$ & $d=3$\\
\hline
\hline
RMSE & \numprint{6.32e-1} &\numprint{7.54e-2}\\
NRMSE &\numprint{4.22e-4} &\numprint{5.03e-5}\\
\hline
\end{tabular}
}
\end{table}

Since usually DEM involve a large amount of data, downsampling is a widely used 
technique to reduce the storage requirements. Our goal is to construct a 
quasi-interpolant spline surface on half of the given data, and then to evaluate 
the produced output on the other half of samples. The obtained results are 
compared with the available \MATLAB{} routines for gridded interpolation: linear, 
cubic, nearest neighbour (N-N), cubic spline (Spline), modified Akima (M-Akima). 
The RMSE and the NRMSE are reported in Table \ref{tab:QIdown}.

\begin{table}[htbp]
\tbl{RMSE and NRMSE results for the downsampling DEM, comparisons with the available \MATLAB{} routines..\label{tab:QIdown}} {
\footnotesize
\centering
\begin{tabular}{|c|c|c|c|c|c|c|}
\hline
\textbf{Linear}  &
\textbf{Cubic}   &
\textbf{Spline}  &
$\mathbf{N}$-N   &
\textbf{M-Akima} & 
\QIBSH{}a &
\QIBSH{}a \\
$d=1$ & $d=3$& $d=3$& $d=0$ & $d=3$& $d=2$& $d=3$\\
\hline
\multicolumn{7}{|c|}{RMSE}\\
\hline
 2.63& 2.10& 2.04& 8.97 &   2.13 & 1.85 &2.05\\
 \hline
\multicolumn{7}{|c|}{NRMSE}\\
\hline
\numprint{1.76e-3} & \numprint{1.41e-3} &  \numprint{1.36e-3}&  \numprint{6.00e-3}&    \numprint{1.42e-3}   &\numprint{1.23e-3} &\numprint{1.37e-3}\\
\hline
\end{tabular}
}
\end{table}

\subsubsection{Curvature Inpainting of Corneal Topographer Data with Missing Regions}\label{test_dati_inpainting}

We present another application  of the tensor product BSH QI  
to a real problem. 
The main goal here is to recover the elevation 
and the radial curvature data of a real eye, processing data 
with missing points.
Sometimes during the topographer acquisition phase
it is impossible to detect data points in some regions:
these are called
inpainting regions and have to be filled in by some numerical techniques.
We apply the TV-H$^{-1}$ inpainting model \cite{burger,schonlieb},
implemented in the~\MATLAB{} function \texttt{bvnegh\_inpainting\_convs}\cite{website:inpainting}.
Applying only this inpainting model to the elevation data was not so succesfull,
hence we  apply it to the radial curvature 
values, combining it with a regularizing phase 
given by the BSH QI interpolant.
We start from the elevation data of the surface
in the regions where they are available, we approximate 
the radial curvature, and then apply digital inpainting 
to the curvature values with missing data.
The curvature is computed along each radial direction approximating the radial partial derivatives 
using the scheme \eqref{eq:finitedifferencescheme}.
After the curvature of the missing regions is recovered using the TV-H$^{-1}$ inpainting, 
the elevations data of the eye are  computed radially, solving a second order ODE. 
The tensor product of the  BSH QI is applied using  its polar coordinates form,
when, during a final step, we regularize the final eye surface. 
Extended discussion, details and numerical results can be found in \cite{inpainting}.

\section{Conclusions}
The  \CPP{} library~\QIBSHPP{} for the approximate solution of several applicative problems provides the implementation of a Hermite type quasi-interpolant operator, with the possibility to approximate the derivatives, with finite difference methods, when they are not available. 
A brief discussion on the theoretical convergence properties of the method is included,
together with some implementation details.
Numerical tests show the convergence properties of the method, even when derivatives are not available and they are approximated.
Computational times and approximation errors make the BSH QI method competitive with other well-known interpolations and 
quasi interpolation methods.
Moreover, the use of the \MATLAB{} C-MEX interfaces for some~\QIBSHPP{} procedures, leads to performance optimization of the code TOM 
in terms of computational time.
The use of the library is also suggested in several applicative fields when high smoothness and high degree splines are required.

\section{Acknowledgments}
The research of Antonella Falini is founded by PON Project AIM 1852414 CUP H95G18000120006 ATT1. 
The authors Antonella Falini and Francesca Mazzia thank the Italian National Group for Scientific Computing (Gruppo Nazionale per il Calcolo Scientifico) for its valuable support under the INDAM-GNCS project CUP\_E55F22000270001.



\bibliographystyle{apalike}
\bibliography{bibliografia}

\end{document}